\documentclass[10pt,reqno]{amsart}
\usepackage{enumerate}
\usepackage{amsfonts}
\usepackage{amsmath}
\usepackage{amsthm}
\usepackage{latexsym}
\usepackage{multicol}
\usepackage{verbatim}
\usepackage{amscd}
\usepackage{amssymb}
\usepackage{hyperref}
\usepackage{amsmath, amscd,wasysym}
\usepackage{graphicx}
\graphicspath{ {images/} }
\usepackage[all,cmtip,arrow]{xy}
\usepackage{pb-diagram}

\advance\textwidth by 1.2in \advance\oddsidemargin by -.6in \advance\evensidemargin by -.6in
\newtheorem*{cor}{Corollary}%[section]
\newtheorem*{lem}{Lemma}
\newtheorem*{prop}{Proposition}

\theoremstyle{definition}

\theoremstyle{definition}
\newtheorem*{thm}{Theorem}

\newcounter{cnt}
 \makeatletter
\def\mydggeometry{\makeatletter\dg@YGRID=1\dg@XGRID=20\unitlength=0.003pt\makeatother}
\makeatother \theoremstyle{remark}

\numberwithin{equation}{section}

  \DeclareMathOperator{\ad}{ad}
\let\bwdg\bigwedge
\def\bigwedge{{\textstyle\bwdg}}

\newcommand{\secref}[1]{Section~\ref{#1}}
\newcommand{\lemref}[1]{Lemma~\ref{#1}}
\newcommand{\propref}[1]{Proposition~\ref{#1}}
\newcommand{\corref}[1]{Corollary~\ref{#1}}

\newcommand{\id}{\operatorname{id}}

\newcommand{\sgn}{\operatorname{sgn}}
\newcommand{\wt}{\operatorname{wt}}

\newcommand{\nc}{\newcommand}
\newcommand{\rnc}{\renewcommand}

\nc{\cal}{\mathcal} \nc{\goth}{\mathfrak} \rnc{\bold}{\mathbf}
\nc{\fk}{\mathfrak}
\newcommand{\supp}{\operatorname{supp}}

\renewcommand{\Bbb}{\mathbb}
\nc\bomega{{\mbox{\boldmath $\omega$}}}\nc\bOmega{{\mbox{\boldmath $\Omega$}}} \nc\bpsi{{\mbox{\boldmath $\Psi$}}}
 \nc\balpha{{\mbox{\boldmath $\alpha$}}} \nc\bll{{\mbox{\boldmath $\ell$}}}
 \nc\bpi{{\mbox{\boldmath $\pi$}}}
\nc\bsigma{{\mbox{\boldmath $\sigma$}}} \nc\bcN{{\mbox{\boldmath $\cal{N}$}}} \nc\bcm{{\mbox{\boldmath $\cal{M}$}}} \nc\bLambda{{\mbox{\boldmath
$\Lambda$}}} \nc\blambda{{\mbox{\boldmath
$\lambda$}}}
\nc\bbeta{{\mbox{\boldmath $\beta$}}}
\nc\bast{{\mbox{\boldmath $\ast$}}}

\nc\btilde{{\mbox{\boldmath ${\bold{\cal{H}}}$}}}

\newcommand{\lie}[1]{\mathfrak{#1}}%\newcommand{\bbtilde}[1]{\mathfrak{#1}}
\makeatletter
\def\section{\def\@secnumfont{\mdseries}\@startsection{section}{1}%
  \z@{.7\linespacing\@plus\linespacing}{.5\linespacing}%
  {\normalfont\scshape\centering}}
\def\subsection{\def\@secnumfont{\bfseries}\@startsection{subsection}{2}%
  {\parindent}{.5\linespacing\@plus.7\linespacing}{-.5em}%
  {\normalfont\bfseries}}
\makeatother

 \nc{\Hom}{\operatorname{Hom}}
  \nc{\td}{\operatorname{\tilde{d}}}\nc{\D}{\operatorname{d}}
  \nc{\mode}{\operatorname{mod}}
\nc{\End}{\operatorname{End}} \nc{\wh}[1]{\widehat{#1}} \nc{\Ext}{\operatorname{Ext}} \nc{\ch}{\text{ch}} \nc{\ev}{\operatorname{ev}}
\nc{\Ob}{\operatorname{Ob}} \nc{\soc}{\operatorname{soc}} \nc{\rad}{\operatorname{rad}} \nc{\head}{\operatorname{head}}
\nc{\A}{\operatorname{\bold{b}}}

 \nc{\Cal}{\cal} \nc{\Xp}[1]{X^+(#1)} \nc{\Xm}[1]{X^-(#1)}
\nc{\on}{\operatorname} \nc{\Z}{{\bold Z}} \nc{\J}{{\cal J}} \nc{\C}{{\bold C}} \nc{\Q}{{\bold Q}} \nc{\R}{{\bold R}}

\nc{\N}{{\Bbb N}} \nc\boa{\bold a} \nc\bob{\bold b} \nc\boc{\bold c} \nc\bod{\bold d} \nc\boe{\bold e} \nc\bof{\bold f} \nc\bog{\bold g}
\nc\boh{\bold h} \nc\boi{\bold i} \nc\boj{\bold j} \nc\bok{\bold k} \nc\bol{\bold l} \nc\bom{\bold m} \nc\bon{\bold n} \nc\boo{\bold o}
\nc\bop{\bold p} \nc\boq{\bold q} \nc\bor{\bold r} \nc\bos{\bold s} \nc\boT{\bold t} \nc\boF{\bold F} \nc\bou{\bold u} \nc\bov{\bold v}
\nc\bow{\bold w} \nc\boz{\bold z} \nc\boy{\bold y} \nc\ba{\bold A} \nc\bb{\bold B} \nc\bc{\bold C} \nc\bd{\bold D} \nc\be{\bold E} \nc\bg{\bold
G} \nc\bh{\bold H} \nc\bi{\bold I} \nc\bj{\bold J} \nc\bk{\bold K} \nc\bl{\bold L} \nc\bm{\bold M} \nc\bn{\bold N} \nc\bo{\bold O} \nc\bp{\bold
P} \nc\bq{\bold Q} \nc\br{\bold R} \nc\bs{\bold S} \nc\bt{\bold T} \nc\bu{\bold U} \nc\bv{\bold V} \nc\bw{\bold W} \nc\bz{\bold Z} \nc\bx{\bold
x} \nc\KR{\bold{KR}} \nc\rk{\bold{rk}} \nc\het{\text{ht }}
\nc\fn{{fin}}  \nc\aff{{aff}}  \nc\tr{{tor}}

\nc{\mpp}{\rotatebox[origin=c]{180}{\pm}}

\nc\toa{\tilde a} \nc\tob{\tilde b} \nc\toc{\tilde c} \nc\tod{\tilde d} \nc\toe{\tilde e} \nc\tof{\tilde f} \nc\tog{\tilde g} \nc\toh{\tilde h}
\nc\toi{\tilde i} \nc\toj{\tilde j} \nc\tok{\tilde k} \nc\tol{\tilde l} \nc\tom{\tilde m}  \nc\ton{\tilde n} \nc\too{\tilde o} \nc\toq{\tilde q}
\nc\tor{\tilde r} \nc\tos{\tilde s} \nc\toT{\tilde t} \nc\tou{\tilde u} \nc\tov{\tilde v} \nc\tow{\tilde w} \nc\toz{\tilde z}
\begin{document}
%\hline
\title{Irreducible Integrable Representations of Toroidal Lie algebras}
\author{Tanusree Khandai}

\address{National Institute of Science Education and Research \\ Khurda 752050\\ India}
\email{tanusree@niser.ac.in}

\iffalse \begin{abstract}
In this paper we classify the irreducible integrable representations with finite-dimensional weight spaces of toroidal Lie algebras on which the centre acts non-trivially. We do so using an approach that was used by  V Chari, G. Fourier and the author in \cite{CFK} to classify the integrable representations with finite-dimensional weight spaces of multi-loop Lie algebras. Using a different approach S.Eswara Rao had classified the irreducible integrable representations with finite-dimensional weight spaces of toroidal Lie algebras in \cite{R2}. %Further we establish a necessary and sufficient condition under which two such irreducible modules are isomorphic.
\end{abstract}\fi

\maketitle

\section{Introduction}
A toroidal Lie algebra $\cal T(\lie g)$ associated to a finite-dimensional complex simple Lie algebra $\lie g_{\fn}$ is the universal central extension of the Lie algebra of polynomial maps from  $(\mathbb{C}^\ast)^k$ to $\lie g_\fn$ where $k$ is a positive integer. From the study of these Lie algebras in \cite{AABGP,RM} it is well known that like the affine Lie algebras, toroidal Lie algebras have a set of real and imaginary roots and one can associate with each real root $\beta$ of $\cal T(\lie g)$ a Lie subalgebra $\lie{sl}_2(\beta)$ which is isomorphic to $\lie {sl}_2(\C)$. A $\cal T(\lie g)$-module  is said to be integrable if it is the direct sum of (possibly infinite) finite-dimensional $\lie{sl}_2(\beta)$-modules for all real roots $\beta$ of $\cal T(\lie g)$. Integrable representations of toroidal Lie algebras and their quotients  have been studied in several papers \cite{CL,L,B,BR,PB,FL,R1,R2,R4,R5}. In \cite{R2} the irreducible integrable $\cal T(\lie g)$-modules having finite-dimensional weight spaces have been classified. In this paper we give an alternative proof of the results in \cite{R2} and also give a parametrization of the isomorphism classes of irreducible $\cal T(\lie g)$-modules with finite-dimensional weight spaces. 

In contrast to an affine Kac-Moody Lie algebra $\lie g_{aff}$ which has a one-dimensional center, the toroidal Lie algebras have a $\Z^k$-graded infinite-dimensional center. This is one of the main sources of difficulty in studying the category $\cal I_\fn^\bast$ of integrable $\cal T(\lie g)$-modules with finite-dimensional weight spaces on which the center acts non-trivially. This problem was sorted in \cite{R2} and it was shown that in each graded component of the center there exists at most one element that acts non-trivially on an irreducible $\cal T(\lie g)$-module in $\cal I_\fn^\bast.$ The proof there used results on representations of Heisenberg Lie algebras from \cite{F}. With an understanding of the integral form of the universal enveloping algebra
 $\cal U(\lie g_{aff})$ of $\lie g_\aff$ from \cite{G} and using the action of certain current algebras on a weight vector of an irreducible $\cal T(\lie g)$-module, in this paper we give an alternate proof of this fact. Besides classifying the simple objects, the methods used here facilitate in the study of the homological properties of the objects in the category  $\cal I_\fn^\bast$ in a way streamlined with \cite{CFK}.
 % using methods that involve basic properties of integrable modules in $\cal I_\fn^{\bast}$ and 

The irreducible integrable $\lie g_\aff$-modules with finite-dimensional weight spaces had been classified in \cite{C1} and their properties studied in a number of papers including \cite{CPloop,CP,CPweyl,CG1,R2,VV}.  Using an approach streamlined with \cite{C1,CFK,CPweyl} we study the simple objects in the category $\cal I_\fn^\bast$  and prove that upto twisting by a one-dimensional $\cal T(\lie g)$-module, every irreducible $\cal T(\lie g)$-module in $\cal I_\fn^\bast$ can be uniquely associated with an orbit for the natural action of $(\C^\ast)^{k-1}$ on the set $\Pi$ of  finitely supported functions from $(\C^\ast)^{k-1}$ to $P_{aff}^+$, the set of dominant integrable weights of the affine Kac-Moody Lie algebra associated to $\lie g_{\fn}$.

The paper is organized as follows. After recalling the structure of the toroidal Lie algebras, in \secref{prelim} we state the results on the representation theory of affine Kac-Moody Lie algebras that play an important role in the classification of the integrable irreducible $\cal T(\lie g)$-modules. We then prove some preliminary results on integrable $\cal T(\lie g)$-modules in \secref{integrable} and finally in \secref{irr.rep.tor} we give an alternative proof of the classification of irreducible $\cal T(\lie g)$-modules in $\cal I_\fn^\bast$. We also establish a necessary and sufficient condition for two irreducible $\cal T(\lie g)$-modules in $\cal I_\fn^\bast$ to be isomorphic.   

\vspace{.15cm}
\noindent{{\textbf{Acknowledgement}}}: I would like to thank S. Viswanath and B. Ravinder for helpful discussions and S. Eswara Rao for reading the earlier drafts carefully and pointing out errors in them.

\section[Preliminaries]{Preliminaries} \label{prelim}
\noindent In this section we fix the notations for the paper and recall the explicit realization of $k$-toroidal Lie algebras from \cite{R2,RM}.

\subsection{}Throughout the paper $\C,\R$, $\Z$ and $\N$ shall denote the field of complex numbers, real numbers, the set of integers and the set of natural numbers, $\Z_+$ shall denote the set of non-negative integers and $\C^\ast$ shall denote the set of non-zero complex numbers.  For a commutative associative algebra $\bold A$,  the set of maximal ideals of $\bold{A}$ shall be denoted by $\max \bold{A}$ and for a Lie algebra $\lie a$ the universal enveloping algebra of $\lie a$ shall be denoted by $\cal U(\lie a)$. For $k\in \N$, a k-tuple of integers $(m_1,\cdots,m_k)$ shall be denoted by $\bom$ and given $\bom\in \Z^k$, $\underline{\bom}$ shall denote the $(k-1)$-tuple of integers $(m_2,\cdots,m_k)$. 

\subsection{}  Let $\fk{g}_{fin}$ be a finite-dimensional simple Lie algebra of rank $n$, $\fk{h}_{fin}$ a Cartan subalgebra of $\fk{g}_{fin}$ and  $R_\fn$ the set of roots of $\fk{g}_\fn$ with respect to $\fk{h}_\fn$.  Let  $\{\alpha_i: 1\leq i\leq n\}$ (respectively $\{\alpha_i^{\vee}: 1\leq i\leq n\}$, $\{\omega_i: 1\leq i\leq n\}$) be a set of simple roots (respectively simple coroots and fundamental weights ) of $\lie g_\fn$ with respect to $\lie h_\fn$, $R_\fn^+$ (respectively $Q_{fin}$, $P_{fin}$) be the corresponding set of positive roots (respectively root lattice and weight lattice) and let $\theta$ (respectively $\theta_s$) be the highest root(respectively highest short root) of $R_\fn^+$ when $\lie g_\fn$ is simply-laced (respectively nonsimply-laced). Let $Q_{fin}^+$ and $P_{fin}^+$ be the $\Z_+$ span of the simple roots and fundamental weights of $(\lie g_\fn,\lie h_\fn)$. Let $\Gamma= P_\fn/Q_\fn$. It is well-known that $\Gamma $ is a finite group whose elements are of the form $\omega_i\mod Q_\fn$ for  $i=1,\cdots,n$. 
\vspace{.15cm} 

 Given $\alpha \in R_\fn^\pm$ let $\lie g_\fn^{\pm\alpha}$ denote the corresponding root space and let $x_\alpha^\pm \in \lie g_\fn^{\pm\alpha}$ and $\alpha^\vee\in \lie h_\fn$ be  fixed elements such that $\alpha^\vee = [x_\alpha^+,x_\alpha^-]$ and $[\alpha^\vee, x_\alpha^\pm] = \pm 2 x_\alpha^\pm$.   
For $\lambda\in P_\fn^+$, let $V(\lambda)$ denote the cyclic $\lie g_\fn$-module generated by a weight vector $v_\lambda$ with defining relations:
$$ x_\alpha^+.v_\lambda =0,\ \forall\ \alpha\in R_\fn^+, \hspace{.5cm} h.v_\lambda =\lambda(h)v_\lambda, \ \forall\ h\in\lie h_\fn,\hspace{.5cm} (x_\alpha^-)^{\lambda(\alpha)+1}.v_\lambda=0, \ \forall\ \alpha\in R^+_\fn.$$ 
It is well known that $V(\lambda)$ is an irreducible finite-dimensional $\lie g_\fn$-module with highest weight $\lambda$ and any irreducible finite-dimensional $\lie g_\fn$-module is isomorphic to $V(\lambda)$ for $\lambda\in P^+_\fn.$ 
Given a non-zero weight vector $u$ in a $\lie g_\fn$-module $V$ be shall denote by $\wt_\fn(u)$ the weight of $u$ with respect to the Cartan subalgebra $\lie h_\fn$ of $\lie g_\fn$. 
\vspace{.15cm}

\subsection{} For a positive integer $k$, let $\C[t_1^{\pm 1},\cdots,t_k^{\pm1}]
$ be the Laurent polynomial ring in $k$ commuting variables $t_1,\cdots,t_k$ and for $\bom=(m_1,\cdots,m_k)\in \Z^k$, let $t^\bom$ (respectively $t^{\underline{\bom}}$) denote the element
$t_1^{m_1}\cdots t_k^{m_k}$ (respectively $t_2^{m_2}\cdots t_k^{m_k}$) in $\C[t_1^{\pm1},\cdots,t_k^{\pm1}]$. Let
$L_k(\fk{g}) = \fk{g}_\fn\otimes \C[t_1^{\pm 1},t_2^{\pm 1},\cdots, t_k^{\pm 1}]$ and $\cal{Z}=\bOmega_k/dL_k$ be the space of K$\ddot{\text{a}}$hler differentials  spanned by the set of vectors $\{ t^\bom K_i,\  ~ \bom\in \Z^k, ~ 1\leq i\leq k\}$ together with the  relation, 
$\sum_{i=1}^k r_i t^{\bor} K_i = 0,$ for $\bor\in \Z^k.$
Let $d_i:({L}_k(\fk{g})\oplus \cal Z)\rightarrow ({L}_k(\fk{g})\oplus \cal Z)$, $1\leq i\leq k$ be the $k$ derivations on $L_k(\fk{g}\oplus \cal Z$ given by: 
\begin{eqnarray}\label{derivation1}
d_i(x\otimes t^\bom)= m_i x\otimes t^\bom, \hspace{.25cm} & \hspace{.25cm}
d_i(t^\bom K_j)=m_i t^\bom K_j & \hspace{.4cm} \forall\ \ 1\leq i,j\leq k,\label{derivation2}
\end{eqnarray} and let $D_k$ be the $\C$ linear span of the derivations $d_1,d_2\cdots,d_k$. The $k$-toroidal Lie algebra associated to a simple Lie algebra $\fk{g}_\fn$ is the vector space $\cal T(\lie g) = L_k(\lie g)\oplus \cal Z\oplus D_k$ on which the Lie bracket is defined by \eqref{derivation1} and the following
relations:
\begin{equation}\label{central}
[x\otimes P, y\otimes Q]=[x,y]\otimes PQ +\overline{Q(dP)} (x|y),\hspace{.35cm}[x\otimes P, \omega ] = 0\hspace{.35cm} [\omega, \omega']=0,
\end{equation} where $x,y\in \fk{g}$, $P,Q\in \C[t_1^{\pm 1},\cdots,t_k^{\pm 1}]$, $\omega, \omega'\in \cal{Z}$ and $\overline{Q(dP)}$ is the residue class of $Q dP$ in $\cal{Z}$. Let $\lie h_{tor}:= \lie h_\fn \oplus \cal Z_0\oplus D_k,$ where $\cal Z_0$ is the subspace of $\cal Z$ spanned by the central elements of degree zero. In order to identify $\lie h_\fn^*$ with a subspace of $\lie h_{tor}^*$, an element $\lambda\in \fk{h}_\fn^*$ is extended to an element of $\lie h_{tor}^*$ by setting $\lambda(c)=0 =\lambda(d_i)=0,$ for all $c\in \cal{Z}_0,  1\leq i\leq k.$
For $1\leq i\leq k$, define $\delta_i \in \lie h_{tor}^*$ by
$ \delta_i|_{\fk{h}_\fn+\cal{Z}_0}=0,$ $\delta_i(d_j) =\delta_{ij},$ for $1\leq j \leq k.$ Given $\alpha\in R_\fn$ and $\bom=(m_1,\cdots,m_n)\in \Z^k$, set $\alpha +\delta_\bom = \alpha +\sum_{i=1}^k m_i\delta_i$ and let  $$R_\tr^{re}= \{\alpha+\delta_\bom : \alpha\in R_\fn, \bom\in \Z^k \}, \hspace{.5cm}  \hspace{.25cm} R_\tr^{im} = \{\delta_\bom =\underset{i=1}{\overset{k}{\sum}} m_i\delta_i : \bom\in\Z^k-\{\bold{0}\} \}. $$ $R_\tr^{re}$ and $R_\tr^{im}$ are respectively the set of real and imaginary roots of $\cal T(\lie g)$ and $R_\tr := R_\tr^{re}\cup R_\tr^{im}$ is the set of all roots of $\cal T(\lie g)$ with respect to $\lie h_{tor}.$ The root vector corresponding to a real root $\alpha+\delta_\bom$ is of the form $x_\alpha\otimes t^\bom$ and the root vectors corresponding to an imaginary root $\delta_\bom$ are of the form $h\otimes t^\bom$ with $h\in \lie h_\fn$. Setting $\alpha_{n+i}:=\delta_i-\theta$, for $ i=1,\cdots,k,$ it can be seen that $\Delta_{tor}=\{\alpha_1,\cdots,\alpha_n, \alpha_{n+1},\cdots,\alpha_{n+k}\}$ forms a simple system for $R_\tr$. %For $\gamma\in R_\tr$, let $\cal T_\gamma$ denote the corresponding root space. It is observed that
%$$ \cal T_{\alpha+\delta_\bom} = t^\bom\otimes \lie g^\alpha_\fn, \hspace{.25cm} \text{for}\ \alpha+\delta_\bom\in R_\tr^{re}, \hspace{.6cm} \text{and} \hspace{.6cm}
%\cal T_{\delta_\bom} = t^\bom\otimes \lie h^\alpha_\fn, \hspace{.25cm} \text{for}\  \delta_\bom\in R_\tr^{im}.$$ 
%Clearly, $\cal{T}_{\alpha+\delta_\bom}$ and $\cal{T}_{\delta_\bom}$ are $\lie{H}$-stable subspaces of $\cal{T}(\lie{g})$ and one has the root space decomposition
%$$\cal{T}(\lie{g}) = \lie {H} \oplus (\underset{\gamma\in R_\tr}{\bigoplus}~ \cal{T}_\gamma) .$$ 

Given $\alpha_\bom =\alpha+\delta_\bom\in R_\tr^{re}$, with $\alpha\in R_\fn^+$ and $\bom\in \Z^k$, let $\alpha_\bom^\vee=\alpha^\vee+\frac{2}{|\alpha|^2}\sum_i m_i K_i.$ With the given Lie bracket operation on $\cal T(\lie g)$ it is easy to check that the subalgebra of $\cal T(\lie g)$
spanned by $\{x_\alpha^+\otimes t^\bom, x_\alpha^-\otimes t^{-\bom}, \alpha_\bom^\vee\}$ is isomorphic to $\lie {sl}_2(\C)$ and we denote it by $\lie {sl}_2(\alpha+\delta_\bom)$.
% $$U(\cal{T}(\lie{g})) = U(\cal{N})\otimes U(\cal{H})\otimes U(\cal{N}^+).$$ gives the corresponding decomposition of the universal enveloping algebra of $\cal{T}(\lie{g})$, where $U(\lie{a})$ denotes the universal enveloping algebra of a Lie algebra $\lie{a}.$

\label{roots}

\vspace{.1cm}
Since $\C[t_1^{\pm 1},\cdots,t_k^{\pm1}]$ is a commutative associate algebra with unity over the complex field $\C$, we see that the following lemma is a special case of \cite[Lemma 2.2]{CFK}.

\begin{lem} Let $\lie g_\fn$ be a simple Lie algebra. Then any ideal of $L_k(\lie g)$, $k\in \N$ is of the form $\lie g_\fn\otimes S$, where $S$ is an ideal of $\C[t_1^{\pm1},\cdots,t_k^{\pm1}]$. Further, $$[ \lie g_\fn\otimes {\C[t_1^{\pm1},\cdots,t_k^{\pm1}]}/{S} ,\ \lie g_\fn\otimes {\C[t_1^{\pm1},\cdots,t_k^{\pm1}]}/{S} ] = \lie g_\fn \otimes {\C[t_1^{\pm1},\cdots,t_k^{\pm1}]}/{S}.$$ 
\end{lem}

\vspace{.15cm}

\subsection{} \label{affine} For $k=1$, the Lie algebra $\cal T(\lie g)$ is called an affine Kac-Moody Lie algebra and we denote it by $\lie g_{aff}$. Explicitly $\lie g_{aff}=\lie g_\fn\otimes \C[t_1^{\pm1}]\oplus \C K_1 \oplus \C d_1.$
Owing to the natural ordering in $\Z $, the set of real and imaginary roots of
$\lie g_{aff}$ can be partitioned as follows:
$$\begin{array}{cc}
R^{ {re}^{\pm}}_{aff}= \{\alpha+m\delta_1 : \alpha\in R_\fn, m\in\Z_\pm \backslash\{0\} \}\cup R_\fn^\pm, \hspace{.35cm} & \hspace{.35cm}
R^{{im}^\pm}_{aff} = \{m\delta_1 : m\in\Z_\pm \backslash\{0\} \}.\end{array}
$$
The set $R_{aff}^+ = R_{aff}^{{re}^+} \cup R_{aff}^{{ im}^+}$ (respectively  $R_{aff}^- = R_{aff}^{{re}^-} \cup R_{aff}^{{ im}^-}$ ) is called the set of positive (respectively negative) roots of $\lie g_{aff}$ and $R_{aff} = R_{aff}^+ \cup R_{aff}^-$ is the set of roots of $\lie g_{aff}$. Denoting the root space of $\lie g_{aff}$ corresponding to a root $\gamma\in R_{aff}$ by $\lie g_{aff}^{\gamma}$, set
$\lie n_{aff}^\pm = \underset{\gamma\in R_{aff}^\pm}{\bigoplus}(\lie g_{aff}^\gamma)$  and $\lie h_{aff}=\lie h_\fn\oplus \C K_1 \oplus \C d_1.$  The
set of simple roots $\Delta_{aff}$ and  coroots $\Delta_{aff}^\vee$ of
$\lie g_{aff}$ are respectively given by
$\Delta_{aff} = \{\alpha_1,\cdots,\alpha_n,\alpha_{n+1}=\delta_1 -\theta\}$, and
$ \Delta_{aff}^\vee = \{\alpha_1^\vee,\cdots,\alpha_n^\vee,\alpha_{n+1}^\vee=K_1-\theta^\vee\}.$ 
Let $Q_{aff}$(respectively $Q_{aff}^\vee$) be the root lattice (respectively
coroot lattice) for $\lie g_{aff}$. Let $\Lambda_i$ ($i=1,\cdots,n,n+1$) be the
fundamental weights of $\lie g_{aff}$, that is,
$\langle \Lambda_i,\alpha_j^\vee\rangle =\delta_{ij},$ for $1\leq j\leq n+1$ and
$\Lambda_i(d_1)=0.$ Then one has the decomposition
$$\lie g_{aff} = \lie n^-_{aff}\oplus \lie h_{aff}\oplus \lie n^+_{aff},  \hspace{.5cm} \text{ and} \hspace{.5cm}\lie h_{aff}^* = \lie h_\fn^* \oplus \C\delta_1 \oplus \C\Lambda_{n+1},$$% where $\Lambda_{n+1}, \delta_1 \in \lie h_{aff}^*$ are such that$\Lambda_{n+1}|_{\lie h^*}= 0 = \Lambda_{n+1}(d_1)$, $\Lambda_{n+1}(K_1)=1$ and
where $\delta_1$ is regarded as an element of $\lie h_{aff}^\ast$ by defining $\delta_1 |_{\lie h_\fn} =0 =\delta_1(K_1)=0$, $\delta_1(d_1)=1$.  Thus an element
$\lambda$ in $\lie h_{aff}^*$ can be uniquely written as
$$\lambda = \lambda(K_1)\Lambda_{n+1} +\lambda\vert_{\lie h_\fn} +
\lambda(d_1)\delta_1,$$ where
$\lambda\vert_{\lie h_\fn}$ is the restriction of $\lambda$ to $\lie h_\fn$.

Let $P_{aff} = \sum_{i=1}^{n+1} \Z\Lambda_i +\C \delta_1,$ (respectively $P_{aff}^+ = \sum_{i=1}^{n+1} \Z_+\Lambda_i +\C \delta_1$) be the set of integral
weights(respectively dominant integral weights) of $\lie g_{aff}$. Let $\succeq$ be the partial order on $P_{aff}$ defined 
by $\lambda\succeq \mu$ if $\lambda,\mu\in P_{aff}$ are such that $\lambda-\mu\in \sum_{i=1}^{n+1} \Z_+ \alpha_{i}$. Given $\lambda,\mu\in P_{aff}$ we shall write $\lambda\succ\mu$ whenever $\lambda\succeq\mu$ but $\lambda\neq \mu$.

\vspace{.15cm}

\subsection{}\label{aff.rep} A $\lie g_{aff}$-module $V$ is said to be integrable if it is $\lie h_{aff}$ diagonalizable and the elements $x_\alpha\otimes t_1^n$, with $\alpha\in R_\fn, n\in \Z$ are locally nilpotent on every $v\in V$. The irreducible integrable $\lie g_{aff}$-modules with finite-dimensional weight spaces were classified in \cite{C1,CPloop,CP}. It was proved that they are either standard modules $X(\Lambda)$, restricted duals of standard modules $X^\ast(\Lambda)$ or loop modules $V(\blambda,\boa,b)$ which are described as follows. 

Given $\Lambda\in P^+_{aff}$, a standard module $X(\Lambda)$ is the unique irreducible $\lie g_{aff}$-module with highest weight $\Lambda$ and highest weight vector $v_\Lambda$ such that $X(\Lambda) = \cal U(\lie g_{aff}).v_\Lambda$. Further $v_\Lambda$ satisfies the relation
$$ \lie n^+_{aff}.v_\Lambda =0, \hspace{1cm} h.v_\Lambda =\lambda(h)v_\Lambda, \
\forall\ h\in\lie h_{aff},\hspace{1cm} (x_{\alpha_i}^-)^{\Lambda(\alpha_i)+1}.v_\Lambda=0, \hspace{.25cm} \forall\ 1\leq i\leq n+1.$$ 
The restricted dual $X^\ast(\Lambda)$ of a standard module $X(\Lambda)$ is a
$\lie g_{aff}$-module generated by a weight vector $v_\Lambda^\ast$ satisfying the relations
$$\lie n_{aff}^-.v_\Lambda^\ast =0, \hspace{1cm} h.v_\Lambda^\ast = -\Lambda(h)v_{\lambda}^\ast, \ \forall\ h\in\lie h_{aff}, \hspace{1cm} (x_{\alpha_i}^+)^{\Lambda(\alpha_i)+1}.v_\Lambda^\ast=0, \hspace{.25cm} \forall\ 1\leq i\leq n+1.$$
Finally, for $r\in \Z_+$, $\blambda =(\lambda_1,\cdots,\lambda_r)\in (P_\fn^+)^r$, $\boa =(a_1,\cdots,a_r)\in (\C^\ast)^r,$ and $b\in \C^\ast$ the loop module $V(\blambda,\boa,b)$ is the the vector space $ V(\lambda_1)\otimes \cdots\otimes V(\lambda_r)\otimes \C[t_1^{\pm1}] $ on which the action of $\lie g_{aff}$ is defined as :
$$ K_1.(v_1\otimes\cdots v_r \otimes f)=0, $$
 $$x\otimes t_1^s.(v_1\otimes\cdots\otimes v_r \otimes f) = \sum_{i=1}^r \ a_i^s  v_1\otimes\cdots\otimes x.v_i\otimes \cdots\otimes v_r \otimes ft_1^s, $$ 
$$d(v_1\otimes\cdots v_r \otimes t^p) = s+p v_1\otimes\cdots v_r \otimes t^p,$$ 
 where $v_i\in V(\lambda_i)$ for $1\leq i\leq r$,\ $f\in \C[t_1^{\pm1}]$,\ $x\in \lie g_\fn$ and $p\in \Z$.
Given an integrable $\lie g_{aff}$-module $V$, let $P_{aff}(V)= \{\eta\in P_{aff}: V_\eta\neq 0\}$, where $V_\eta = \{v\in V: hv = \eta(h)v, \ \text{for}\ h\in \lie h_{aff}\}$ for $\eta\in P_{aff}$. For $v\in V_{\eta}$ we write $\wt_{aff}(v) = \eta.$ With this notation we state the following results. Parts (i) and (ii) of the following proposition have been proved in \cite{CPloop,C1} and parts (iii) and (iv) have been proved in \cite{RedR,CG1}.

\vspace{.15cm}

\begin{prop} 
\item[(i).]  Let $r\in \N$. Given $\blambda \in (P_\fn^+)^r, \boa\in (\C^\ast)^r$ and $b\in \C$, the loop module $V(\blambda,\boa,b)$ is completely reducible as a $\lie g_{aff}$-module.
  %Given $\upsilon\in \Upsilon$, there exists $r_\upsilon\in \Z_+$ such that  $L^s(V_\upsilon)$ is an irreducible $\lie g_{aff}$ module for  $0\leq s< r_\upsilon$ and 
   %$$V^g_\upsilon = \bigoplus_{s=0}^{r_\upsilon-1} L^s(V_\upsilon).$$
  % \item[2.] For every $a\in \C^\times$ and $\mu\in P_{fin}^+$,  $V^g_{\upsilon_{a,\mu}} = V_{\upsilon_{a,\mu}}\otimes \C[t_1^\pm]$ is an irreducible $\lie g_{aff}$-module.
\item[(ii).]  Let $V$ be an irreducible representation of $\lie g_{aff}$ having finite-dimensional weight spaces. If $m\in \Z$ is such that $K_1.v = m v$, for all $v\in V$, then \begin{itemize}
\item[a.] for $m>0$, (respectively $m<0$)  $V$ is isomorphic to $X(\Lambda)$ (respectively $X^\ast(\Lambda)$) for some $\Lambda\in P_{aff}^+$.
\item[b.] for $m=0$, $V$ is isomorphic to an irreducible component of a loop module $V(\blambda,\boa,b)$ for some $\blambda \in (P_\fn^+)^r, \boa\in (\C^\ast)^r, b\in \C$ and $r\in \Z$.\end{itemize}
\item[(iii).]If all eigenvalues for the action of the central element $K_1$ on an integrable $\lie g_{aff}$-module $V$, are non-zero then $V$ is completely reducible as $\lie g_{aff}$-module.
\item[(iv).] Let $\lambda\in P_{aff}^+$ be of the form $\lambda = \lambda(K_1)\Lambda_{n+1} + \lambda|_{\lie h_\fn}+\lambda(d_1)\delta_1$. If $\varpi_\lambda$ is the  unique minimal element in $P_{fin}^+$ such that $\lambda|_{\lie h_\fn} \equiv \varpi_\lambda\mod Q_{fin}$ then $\varpi_{\lambda,r} = \lambda(K_1)\Lambda_{n+1}+\varpi_\lambda+r\delta_1\in P_{aff}(X(\lambda))$ for all $r\in \C$ such that $\lambda(d)-r\in \Z_+$. In particular, if
  $\lambda|_{\lie h_\fn}\in P_{fin}^+\cap Q_{fin}^+$ and $\lambda|_{\lie h_\fn}\neq 0$, then $\lambda(K_1)\Lambda_{n+1}+\beta+r\delta_1 \in P_{aff}(\lambda)$ for all $r\in \C$ such that $\lambda(d_1)-r\in \Z_+$, where $\beta=\theta$ if $\lie g$ is simply-laced and $\beta=\theta_s$ otherwise.
 \end{prop}  

\vspace{.15cm}

\subsection{}\label{Garland}Let $u$ be an indeterminate. For $\alpha \in R_\fn$ define a power series $\bold{p}_{\alpha}(u)$ in $u$ with coefficients in $\cal U(\lie h_\fn\otimes \C[t^{\pm1}])$ :
$$\bold{p}_{\alpha}(u)= \exp \left(-\underset{r=1}{\overset{\infty}{\sum}}\frac{\alpha^\vee\otimes t^r}{r}u^r\right).$$ For $s\in \Z_+$, let $p_\alpha^s$ be the coefficient of $u^s$ in $\bold{p}_{\alpha}(u).$ The following was proved in \cite[Lemma 7.5]{G}.

\begin{lem} Let $\alpha\in R_\fn^+$. Then for $r \geq 1$ we have
\begin{equation}
(x_\alpha^+\otimes t)^r(x_\alpha^-\otimes 1)^{r+1}-\underset{s=0}{\overset{r}{\sum}}(x_\alpha^-\otimes t^{r-s})p_{\alpha}^s \in \cal U(\lie g_{aff})(\lie n_{aff}^+),
\end{equation}
\begin{equation}
(x_\alpha^+\otimes t)^{r+1}(x_\alpha^-\otimes 1)^{r+1}-p_{\alpha}^{r+1} \in \cal U(\lie g_{aff})(\lie n_{aff}^+).
\end{equation} 
 \end{lem}

\vspace{.15cm}

\section{Integrable Representations of the toroidal Lie algebra}
\label{integrable}
\subsection{} Let $\cal I$ be the category whose objects are integrable $\cal{T}(\lie g)$-modules and morphisms
$$\Hom_{\cal I}(V,V^\prime) = \Hom_{\cal T(\lie g)}(V, V^\prime), \hspace{.35cm} V,V^\prime \in \cal I.$$
For an integrable representation   $\psi:{\cal T}(\lie g)\rightarrow End(V)$  of ${\cal T}(\lie g)$, let $V^\psi$ denote the corresponding $\cal{T}(\lie g)$-module in $\cal I$. Given $\beta = \alpha+\delta_\bom\in R_\tr^{re} $ with $\alpha\in R_\fn$, define an operator $r_{\beta}^\psi$ on $V$ as follows:
$$
r_\beta^\psi = \exp(\psi(x_{\alpha}^+ \otimes t^\bom))
\exp(\psi(-x_{\alpha}^-\otimes t^{-\bom}))\exp(\psi(x_{\alpha}^+\otimes t^\bom)).
$$
Since $V^\psi$ is an integrable $\cal T(\lie g)$-module, the operator $r_\beta^\psi$ is well-defined for all $\beta\in R_\tr^{re}$. Let $$P(V^\psi) = \{\lambda\in \lie h_\tr^\ast : V^\psi_{\lambda}\ne 0\}, \hspace{1cm}\text{where }\ \ V^\psi_{\lambda} =\{v\in V^\psi: h.v=\lambda(h)v, \text{for}\ h \in \lie h_\tr\}.$$ Let $W^\psi_\tr = \langle r_\beta^\psi: \beta\in R_\tr^{re} \rangle, $ be the group generated by the operators $r^\psi_\beta$ for $\beta\in R_\tr^{re}$.  By \cite[Lemma 3.8, \S6.5]{Kac}, $r_\beta^\psi(\lambda)= \lambda - \langle\lambda,\beta^\vee\rangle\beta,$ for $\lambda\in P(V^\psi)$. Using the representation theory of $\lie{sl}_2(\C)$ the following is standard in $\cal I.$

\label{sl2rep}

\vspace{.15cm}
\begin{lem}Let $V^\psi$ be a $\cal T(\lie g)$-module in $\cal I$ and let $\lambda\in P(V^\psi)$. Then,
\begin{itemize}
\item[i.]$\langle \lambda, \alpha^\vee\rangle\in \Z$, for $\alpha\in R_\tr^{re}.$
\item[ii.] $w\lambda \in P(V^\psi)$ and $\dim V_\lambda^\pi = \dim V_{w\lambda}^\pi$ for all $\lambda\in P(V^\psi)$, $w\in W^\psi_\tr$.
% \item[iii.]If $\lambda+\alpha\notin P(V^\psi)$ (respectively, $\lambda-\alpha\notin P(V^\psi)$) then $\langle \lambda,\alpha^\vee\rangle\geq 0$ (respectively, $\langle \lambda,\alpha^\vee\rangle\leq 0$), for $\lambda\in P(V^\psi), \alpha\in R_\tr^{re}$.
%\item[iv.] Given $\lambda\in P(V^\psi)$ and $\gamma\in R_\tr^{re}$ such that $\langle \lambda,\gamma^\vee\rangle> 0$ (respectively, $\langle \lambda,\gamma^\vee\rangle< 0$), $\lambda-\gamma\in P(V^\psi)$ (respectively, $\lambda+\gamma\in P(V^\psi)$).
\item[iii.] Let $\alpha\in R_\fn^+$ and $\beta=\alpha+m_i\delta_i \in R_\tr^{re}.$ Given a $\cal T(\lie g)$-module $V^\psi$ in $\cal I$ and $\lambda\in P(V^\psi)$, we have $$r_\alpha^\psi r_\beta^\psi(\lambda) = \lambda + \frac{2}{(\alpha|\alpha)}( m_i\langle\lambda, K_i\rangle) \alpha - (\langle\lambda,\alpha^\vee\rangle+\frac{2}{(\alpha|\alpha)} m_i\langle\lambda, K_i\rangle)\delta_i.$$ \end{itemize}
\end{lem}

The following is an easy corollary of \lemref{sl2rep}(iii).
\begin{cor} Let $V^\psi$ be an integrable $\cal T(\lie g)$-module. Suppose $\lambda+\underset{i=1}{\overset{k}{\sum}} r_i\delta_i\in P(V^\psi)$ satisfies the condition \begin{equation}\label{m.00}\langle \lambda,\alpha_{n+1}^\vee\rangle = m, \hspace{1cm}\text{and}\hspace{1cm} \langle \lambda,K_j\rangle=0, \hspace{.25cm} \text{for}\ j=2,\cdots,k.\end{equation} Then there exists $\underline{\bom}=(m_2,\cdots,m_k)\in \Z^{k-1}$ with $0\leq m_i< m$ for $2\leq i\leq k$ such that $$\lambda +r_1\delta_1 +\underset{i=2}{\overset{k}{\sum}}  m_i\delta_i \in P(V^\psi).$$
\end{cor}
\proof Set $\beta_j=\delta_j + \delta_1 -\theta = \delta_j+\alpha_{n+1}$ for $j=2,\cdots,k$. Then by \ref{roots}, $$\beta_j^\vee = -\theta^\vee +\frac{2}{(\theta|\theta)}(K_j+K_1) = +\alpha_{n+1}^\vee +\frac{2}{(\theta|\theta)}K_j, \hspace{.5cm} 2\leq j\leq k.$$  By definition $\alpha_{n+1}, \beta_2, \beta_3,\cdots, \beta_k \in R_\tr^{re}$, hence $r_{\alpha_{n+1}}^\psi, r_{\beta_2}^\psi,\cdots, r_{\beta_k}^{\psi}$ are well-defined operators in $W^\psi_\tr$. Thus given $\lambda+\underset{i=1}{\overset{k}{\sum}} r_i\delta_i\in P(V^\psi)$, by \lemref{sl2rep}(ii), $r_{\alpha_{n+1}}^\psi r_{\beta_j}^\psi(\lambda+\underset{i=1}{\overset{k}{\sum}} r_i\delta_i) \in P(V^\psi)$ for $2\le j\leq k$. Now note that
$$
r^\psi_{\alpha_{n+1}}r^\psi_{\beta_j}(\lambda+\delta_\bor) = r^\psi_{\alpha_{n+1}}(\lambda+\delta_\bor -\langle \lambda+\delta_\bor, \beta_j^\vee\rangle \beta_j)= r^\psi_{\alpha_{n+1}}(\lambda+\delta_\bor -\langle \lambda+\delta_\bor, \frac{2}{(\theta|\theta)} K_j+\alpha_{n+1}^\vee \rangle (\delta_j +\alpha_{n+1})$$
$$=  r^\psi_{\alpha_{n+1}}(\lambda+\delta_\bor -\langle \lambda+\delta_\bor, \alpha_{n+1}^\vee \rangle \alpha_{n+1}) + \langle \lambda +\delta_\bor, \frac{2}{(\theta|\theta)}K_j \rangle \alpha_{n+1} -\langle \lambda+\delta_\bor, \frac{2}{(\theta|\theta)} K_j+\alpha_{n+1}^\vee \rangle \delta_j $$
$$ = (\lambda+\delta_\bor) +\langle \lambda +\delta_\bor, \frac{2}{(\theta|\theta)}K_j \rangle (\alpha_{n+1}-\delta_j)-\langle \lambda +\delta_\bor, \alpha_{n+1}^\vee\rangle \delta_j. 
$$ Since  $\langle \lambda, \alpha_{n+1}^\vee\rangle =m$ and $\langle \lambda, K_j\rangle =0$ for $2\leq j\leq k$, it follows that $\lambda+\delta_\bor-m\delta_j \in P(V^\psi)$. Repeating the argument it is easy to see that if $\lambda+\delta_\bor\in V^\psi$ satisfies the condition \eqref{m.00} and $\bor\in \Z^k$ is such that $r_j\geq 0$ for $2\leq j\leq k$, then there exists $w\in W^\psi_\tr$ such that $w(\lambda+\delta_\bor) = \lambda +r_1\delta_1 +\sum^k_{i=2} m_i\delta_i$ with $0\leq m_i<m$ for $2\leq i\leq k.$ If $\bor\in \Z^k$ is such that $r_j\leq 0$ for some $2\leq i\leq k$ then the same result can be obtained by setting $\gamma_j=\delta_j-\alpha_{n+1}$ and using the operator $r_{\gamma_j}$ in place of $r_{\beta_j}$.\endproof

\vspace{.15cm}

Assuming that we understand that the elements of $\cal{T}(\lie g)$ act on an object $V$ of $\cal I$  via a Lie algebra homomorphism $\psi: \cal{T}(\lie g)\rightarrow End(V)$, we shall henceforth drop the superscript $\psi$ when referring to $V^\psi$, $W_\tr^\psi$ etc.

\vspace{.35cm}

\subsection{} Let $\cal Z_0$ be the $\C$-span of the zero degree central elements of $\cal T(\lie g)$. Given a $\cal T(\lie g)$-module $V$ and $\Lambda\in P(V)$, the restriction $\Lambda|_{\cal Z_0}$ is a map from $\cal Z_0$ to $\C$. In particular when $V$ is integrable,  $\Lambda|_{\cal Z_0}\subset\Z$. For $\bon=(n_1,\cdots,n_k)\in \Z^k$ let
$P_{\bon}(V)=\{\lambda\in P(V): \lambda(K_i)=n_i, \text{for}\ 1\leq i\leq k\}$ and let
$$V^{(\bon)} = \underset{\lambda\in P_{\bon}(V)}{\bigoplus} V_\lambda.$$
As $\cal{Z}_0$ commutes with $\cal {T}(\lie g)$, $V^{(\bon)}$ is a $\cal{T}(\lie g)$-module for each $\bon\in \Z^k$ and any $V\in \Ob \cal I$ can be decomposed as
follows : $$V= \underset{\bon\in \Z^k}{\oplus} \big(\underset{\lambda\in P_\bon(V)}{\bigoplus} V_\lambda\big) = \underset{\bon\in \Z^k}{\oplus} V^{(\bon)}.$$ 
Note that for all $\lambda\in P_{\bon}(V)$, the restriction $\lambda|_{\cal {Z}_0}$ is a linear functional. Hence by a change of basis for $\cal Z_0$ one can always assume that at most one zero degree central element acts non-trivially on $V$. Since in a toroidal Lie algebra the choice of basis for $\cal Z_0$ is dependent on the choice of generators of the coordinate ring $\C[t_1^{\pm 1},\cdots,t_k^{\pm 1}]$,  the change of basis matrix $\bold{B} = (b_{ij})_{1\leq i,j\leq k}$ for $\cal Z_0$ induces a homomorphism $\tilde{\bold{B}}: L_k(\lie g)\rightarrow L_k(\lie g)$ given by : $$\tilde{\bold{B}}(x\otimes t^\bom) = x\otimes t^{\bom\bold{B}}.$$
Let $\{\boe_i: 1\leq i\leq k\}$ be the standard basis of $\mathbb{R}^k.$ Then setting $a_i =t^{\boe_i\bold{B}}$ for $1\leq i\leq k$, it is easy to see that $\tilde{\bold{B}}: L_k(\lie g)\rightarrow \lie g_\fn\otimes \C[a_1^{\pm 1},\cdots,a_k^{\pm 1}]$ is an isomorphism of multiloop Lie algebras whenever $\bold{B}$ is an integer matrix with determinant $\pm 1$. As noted in \cite{R1}, the isomorphism $\tilde{\bold{B}}$ can be extended to an isomorphism of toroidal Lie algebras $$\tilde{\A}: \cal T(\lie g)\rightarrow \lie g_\fn\otimes \C[a_1^{\pm 1},\cdots,a_k^{\pm 1}]\oplus \tilde{\cal Z}\oplus \tilde{D_k},$$ where $\tilde{\cal Z}$ is the $\C$ span of  $\{a^\bom K_i': \sum_{i=1}^k r_i a^\bor K_i'=0, \bom,\bor\in \Z^k, 1\leq i\leq k\},$ with $K_i' = \sum_{j=1}^k b_{ij}K_i$ for $1\leq i\leq k$ and $\tilde{D_k}$ is the $\C$-span of the derivations $\{\tilde{d_i}= a_i\frac{d}{da_i}: 1\leq i\leq k\}$. Defining the $\cal T(\lie g)$ action on $V^{(\bon)}\in \Ob \cal I$ by
$$X.v = \tilde{\A}(X).v, \hspace{.5cm} \forall\ v\in V,$$ it is then easy to see that with respect to the toroidal Lie algebra $\tilde{\A}(\cal T(\lie g)),$ the module $V^{(\bon)}$ is of the form $V^{(m\boe_1)}$ for some $m\in \Z.$ The following standard result from group theory shows that this is true in general.i.e., upto an isomorphism any indecomposable object in $\cal I$ on which the center acts non-trivially is of the form $V=V^{(m\boe_1)}$ for some $m\in \Z$.

\begin{lem} Let $\{n_1,\cdots,n_k\}$ be a set of  integers with $\gcd(n_1,\cdots,n_k)\neq 0.$ Then there exists a $k\times k$ integer matrix $\bold{B}=(b_{ij})$ with determinant $\pm 1$ such that
  \begin{equation}\label{three-one}
\sum_{j=1}^k b_{1j}n_j = \gcd(n_1,\cdots,n_k), \hspace{.5cm}\text{and}\hspace{.35cm}  \hspace{.5cm} \sum_{j=1}^k b_{ij} n_j = 0, ~ \forall\ i\geq 2,
\end{equation}
  \end{lem}

\label{zero.degree.center}

\vspace{.1cm}

\subsection{}For $\bom\in \Z^k$ let  $\cal I^{(\bom)}$ be the full subcategory of $\cal I$ whose objects are $\cal T(\lie g)$-modules on which the zero degree central element $K_i$,  acts by the integer $m_i$ for $1\leq i\leq k$. The following is an immediate consequence of \lemref{zero.degree.center} .

 \label{three-five}

 \begin{lem} Let $V$ be an integrable $\cal T(\lie g)$-module. Then $$V= \underset{\bom\in \Z^k}{\bigoplus} V^{(\bom)}, \hspace{1cm} \text{where}\ V^{(\bom)}= \underset{\lambda \in P_\bom(V)}{\oplus} V_\lambda.$$ Upto an isomorphism the component $V^{(\bom)}$ of $V$ is of the form $V^{(m\boe_1)}$ where $m=\gcd(m_1,\cdots,m_k).$ Further, $\Ext_{\cal I}^1(V, U)=0$ for $V\in \Ob \cal I^{(\bom)},\ U\in \Ob \cal I^{(\bon)}$ with $\bom,\bon\in Z^k$ whenever $\bom\neq \bon$. In particular,
   $$\cal I = \underset{\bom\in \Z^k}{\bigoplus} \cal I^{(\bom)},$$ and  given 
$   \bom=(m_1,\cdots,m_k)\in \Z^k$    the category $\cal I^{(\bom)}$ is equivalent to $\cal I^{(m\boe_1)}$ where $m=\gcd(m_1,\cdots,m_k)$.\end{lem} 
\noindent Without loss of generality we thus restrict ourselves to the study of the full subcategory $\cal I^{(m\boe_1)}$ of $\cal I$.% Henceforth by an integrable $\cal T(\lie g)$-module we shall refer to an object in $\cal I^{(m\boe_1)}.$

 \subsection{}Let $\cal I_{fin}$ ( respectively $\cal I_\fn^{(m\boe_1)}$) be the full subcategory of $\cal I$ ( respectively $\cal I^{(m\boe_1)}$) consisting
 of integrable $\cal{T}(\lie g)$-modules with finite dimensional weight spaces.
% and for $m\in \Z$, let $\cal I_{fin}^{(m)}$ be the full subcategory of  $\ca I^{(m\boe_1)}$ consisting of integrable $\cal{T}(\lie g)$-modules with
%finite dimensional weight spaces.

 \begin{prop}\label{highest.wt} Let  $V$ be an integrable $\cal T(\lie g)$-module in  $\cal I_{fin}^{(m\boe_1)}$, where $m>0$. Let $\lie n_{aff}^+$ be the positive root space of the affine Lie subalgebra $\lie g_{aff} = \lie g_\fn\otimes \C[t_1^{\pm 1}]\oplus \C K_1\oplus \C d_1$ of $\cal T(\lie g)$. Then
$$V^+_{aff}=\{v\in V_\lambda : \lie n^+_{aff}\otimes \C[t_2^{\pm 1},\cdots,t_k^{\pm 1}].v=0\}$$ is a non-empty subset of $V$.
\end{prop}
\proof  For a $(k-1)$ tuple $\underline{\bor}=(r_2,r_3,\cdots,r_k)\in \Z^{k-1}$, let
$$V[\underline{\bor}] =\{v\in V_\lambda:d_i v= r_i v,\ 2\leq i\leq k,\ \lambda\in P(V)\}.$$ Clearly $V[\underline{\bor}]$ is an integrable $\lie g_{aff}$-module with finite-dimensional weight spaces. Using \ref{roots}, we can write
$$V[\underline{\bor}] = \underset{\gamma\in P_\fn/Q_\fn}{\bigoplus} V^\gamma[\underline{\bor}], \hspace{1.35cm} \text{where} \hspace{.15cm} V^\gamma[\underline{\bor}] = \underset{\varpi_\gamma \equiv \lambda|_{\lie h_\fn}\mod Q_\fn}{\bigoplus} (V_{\lambda}\cap V[\underline{\bor}]).$$
Since $m>0$, by \propref{aff.rep}(iii), $V^\gamma[\underline{\bor}]$ is completely reducible as a $\lie g_{aff}$-module for each $\gamma\in P_\fn/Q_\fn$. Hence there exists a (possibly infinite) $\lie g_{aff}$-module filtration $\cdots \subset V_1\subset V_{0} = V^\gamma[\underline{\bor}]$ of  $V^\gamma[\underline{\bor}]$ such that the successive quotients are isomorphic to the $\lie g_{aff}$-module summands of $V^\gamma[\underline{\bor}]$.  %. Since $V[\bor^1]\{\gamma\}$ is completely reducible a filtration can be chosen such that
Without loss of generality we may assume that $V_j/V_{j+1}\cong_{\lie g_{aff}} X(\lambda_j)$ with $\lambda_j\in P^+_{aff}$. % and $\lambda_i(d_1)\geq \lambda_j(d_1)$ for $i>j$. Then 
Setting $\tilde\lambda_j= \lambda_j+\sum_{i=2}^k r_i\delta_i$ we see that there exists $v_j\in (V_j/V_{j+1})_{\tilde{\lambda}_j}$, such that $\cal{U}(\lie g_{aff})v_j$ is isomorphic as a $\lie g_{aff}$-module to $X(\lambda_j)$. Since $\lambda_j|_{\lie h_\fn}\equiv\varpi_\gamma\mod Q_\fn$, by \propref{aff.rep}(iv), $m\Lambda_{n+1}+\varpi_\gamma+\lambda_j(d_1)\delta_1$ is a weight of $X(\lambda_j)$. Let $r_1^j = \underset{1\leq s\leq j}{\min} \lambda_s(d_1).$ Then by \lemref{aff.rep}(iv)
$$\Lambda_{\gamma, \underline{\bor}}^j =m\Lambda_{n+1}+\varpi_\gamma+r_1^j\delta_1+\sum_{i=2}^k r_i\delta_i \in P(V_s/V_{s+1})$$ for all $s\in \Z_+$ such that $s\leq j$, whenever $\lambda_j|_{\lie h_\fn}\not\in Q_\fn^+$ and $$\Lambda_{0,\underline{\bor}}^j=m\Lambda_{n+1}+\lambda_j(d_1)\delta_1+\sum_{i=2}^k r_i\delta_i \in P(V_s/V_{s+1})$$ for all $s\in \Z_+$ such that $s\leq j$, whenever $\lambda_j|_{\lie h_\fn}\in P^+_\fn\cap Q^+_\fn$. Let $x_n = \dim V_{\Lambda^n_{\gamma,\underline{\bor}}}$(respectively $V_{\Lambda^j_{0,\underline{\bor}}}$). From above argument it follows that $(x_n)_{n\in \N}$ is a sequence of integers such that $x_j\geq j$ for all $j\in \N$. In particular  
for every $K\in \N$, $x_n> K$ for all $n\geq K+1$. Thus if  for some $\underline{\bor}\in \Z^{k-1}$, $V^\gamma[\underline{\bor}]$ would have an infinite filtration, %n setting $x_n = \dim V_{\Lambda^n_{\gamma,\underline{\bor}}}$(respectively $V_{\Lambda^j_{0,\underline{\bor}}}$) we would get a sequence of integers $(x_n)_{n\in \N}$ such that for each $j\in \N$,  %we would have 
% $$x_j = \dim V_{\Lambda_{\gamma,\underline{\bor}}^j} \geq j \hspace{.5cm}(\text{respectively}\ \dim V_{\Lambda^j_{0,\underline{\bor}}} \geq j ).$$ Clearly the sequence $(x_n = dim V_{\Lambda_{\gamma,\underline{\bor}}^n})$ tends to infinity. %Hence  %Thus if $V[\underline{\bor}]$ admits an infinite filtration 
then the fact that $\underset{n\rightarrow \infty}{\lim} x_n=\infty$ would contradict the finite-dimensionality of the weight spaces of $V$. % $V_{\Lambda_{\gamma,\underline{\bor}}}$  or $V_{\Lambda_{0,\underline{\bor}}}$ (as the case may be) would be contradicted. 
Hence $V[\underline{\bor}]$ admits a finite $\lie g_{aff}$-module filtration for every $\underline{\bor}=(r_2,\cdots,r_k)\in \Z^{k-1}$. Consequently for each $\gamma\in P_\fn/Q_\fn $ and $\underline{\bor} =(r_2,\cdots,r_k)\in \Z^{k-1}$ there exists $\lambda_{\gamma,\underline{\bor}}\in P_{aff}^+$ such that $X(\lambda_{\gamma, \underline{\bor}})$ is a summand of $V^\gamma[\underline{\bor}]$ but $X(\lambda_{\gamma,\underline{\bor}}+\beta)$ is not a summand of $V^\gamma[\underline{\bor}]$ for any $\beta\in R_{aff}^+$.

Let $$\Z^{k-1}_m=\{\underline{\bom}=(m_2,\cdots,m_k)\in \Z^{k-1}: 0\leq m_i<m,\ 2\leq i\leq k\}.$$ Since $\Z^{k-1}_m$ is a finite subset of $\Z^{k-1}$, given $\gamma\in P_\fn/Q_\fn$ there exists $\lambda_\gamma\in P_{aff}^+$ such that $X(\lambda_\gamma)$ is a summand of $V^\gamma[\underline{\bom}]$ for some $\underline{\bom}\in \Z^{k-1}_m$ and $\lambda_\gamma\succeq \lambda_{\gamma,\underline{\bom}}$ for all $\underline{\bom}\in \Z^{k-1}_m.$ On the other hand by \corref{sl2rep}, for $\lambda\in P_{aff}^+$, $\lambda+ \delta_{\underline{\bor}} \in P(V^\gamma)$ if and only if $\lambda +\delta_{\underline{\bom}}\in P(V^\gamma)$ for some $\underline{\bom}\in \Z^{k-1}_m.$ Hence it follows that $\lambda_\gamma\succeq \lambda_{\gamma,\underline{\bor}}$ for all $\underline{\bor}\in \Z^{k-1}$ implying that there exists a non-zero vector $v\in V_{\lambda_\gamma +\sum_{i=2}^k m_i\delta_i},$ with  $\underline{\bom}\in \Z^{k-1}_m$ such that  $\lie n_{aff}^+\otimes \C[t_2^{\pm 1},\cdots,t_k^{\pm 1}].v=0.$ This completes the proof of the proposition.
\endproof

\section{Irreducible objects in $\cal I_{fin}^{(m\boe_1)}, m>0$}
\label{irr.rep.tor}
In this section we give a parametrization for the irreducible $\cal T(\lie g)$-modules in $\cal I_{fin}^{(m\boe_1)}.$

\subsection{} The following proposition, which was first proved in \cite{R2}, plays an important role in determining the simple objects in $\cal I_\fn^{(m\boe_1)}, m>0$. We give here an alternative proof using an understanding of the integral forms in the universal enveloping algebra $\cal U(\cal T(\lie g))$ of $\cal T(\lie g)$.

\begin{prop} Let $V$ be an irreducible $\cal T(\lie g)$-module in $\cal I_{fin}^{(m\boe_1)}$, $m>0$.   Then  $K_jt^\bor$ acts trivially on $V$ for all $2\leq j\leq k$ and all monomials  $t^\bor \in \C[t_1^\pm,t_2^\pm,\cdots,t_k^\pm]$.
%\item[(i).]  
%\item[(ii).] $K_1t^\bor$ acts trivially on $V_{aff}^+$ whenever $\bor=(r_1,\cdots,r_k)\in \Z^k$ is such that $r_1>0$.
%\item[(ii).] If $V$ is irreducible then in addition to $K_jt^\bor$, $j\geq 2$, the central elements  %the graded central elements    $K_jt^\bor$ acts trivially on $V$ for all $2\leq j\leq k$ and all monomials  $t^\bor \in \C[t_1,t_2^\pm,\cdots,t_k^\pm]$. 
Further, $K_1t^\bor$ acts trivially on $V$ whenever $\bor=(r_1,\cdots,r_k)\in \Z^k$ is such that $r_1\neq 0.$
\end{prop}
\proof Given $1\leq j\leq k$, and  $\bor\in \Z^k-\{\bold{0}\}$ let $\cal N_{j,\bor} = \{u\in V: K_jt^\bor.u=0\}$. Clearly $\cal N_{j,\bor}$ is a submodule of the irreducible $\cal T(\lie g)$-module $V$. Hence to prove the proposition it suffices to show that $\cal N_{j,\bor}$ is non-empty for each $\bor\in \Z^k$ whenever $2\leq j\leq k$ and $\cal N_{1,\bor}$ is non-empty whenever $r_1\neq 0$.

We prove the proposition in 3 steps. Let $\bor=(r_1,\cdots,r_k)\in \Z^k.$ In step 1 we show that $\cal N_{j,\bor}\neq 0$ for all $1\leq j\leq k$ whenever $r_1>0$. In step 2 we show that $\cal N_{j,\bor} \neq 0$ for $2\leq j\leq k$ whenever $r_j\neq 0$ and in step 3 we take care of the remaining cases. That is we show that $\cal N_{j,\bor} \neq 0$ for $2\leq j\leq k$ whenever $r_j=0$ and as a consequence $\cal N_{1,\bor}\neq 0$ whenever $r_1<0$.

\noindent {\bf{Step 1}}. By \propref{highest.wt} the subspace $V_{aff}^+$ of $V$ is non-empty and for any $v\in V_{aff}^+$ we have
\begin{equation}\label{n.aff.+}\begin{array}{l} h\otimes t^{\bop}.v =0, \hspace{.5cm} \forall\ \bop =(p_1,p_2,\cdots,p_k) \in \Z^k\  \text{with}\ p_1>0,\ h\in \lie h_\fn,\\
x_\alpha^+\otimes t_i^s.v =0 \hspace{.5cm} \forall\ s\in \Z,\ \alpha\in R_\fn^+,\
2\leq i\leq k,\\
 x_\alpha^-\otimes t^\bor.v =0\hspace{.5cm} \forall\ \bop =(p_1,p_2,\cdots,p_k) \in \Z^k\  \text{with}\ p_1>0, \ \alpha\in R_\fn^+.
 \end{array} 
 \end{equation} Hence given $\alpha\in R_\fn^+$ for all $2\leq i\leq k$ we have
$$0=(x_\alpha^+\otimes t_i^s)( x_\alpha^-\otimes t^\bor)v = (x_\alpha^- \otimes t^\bor)( x_\alpha^+\otimes t_i^s)v + (\alpha^\vee \otimes t^\bor t_i^{s}+K_it^\bor t_i^s).v,$$ whenever $\bor=(r_1,\cdots,r_k)\in \Z^k$ is such that $r_1>0$.
Using the relations \eqref{n.aff.+} it thus follows that $K_i t^\bos.v =0$ for all $2\leq i\leq k$ whenever $\bos=(s_1,\cdots,s_k)$ is such that $s_1>0$, i.e., $\cal N_{i,\bos}\neq 0$ for $2\leq i\leq k$ whenever $s_1>0$. As 
$\sum_{j=1}^k s_j K_jt^{\bos} =0$ it follows that $\cal N_{1,\bos}\neq 0$ whenever $s_1>0$.

\noindent {\bf{Step 2}}. For a contradiction assume that there exists $2\leq j\leq k$ and $\bor\in \Z^k$ with $r_1 =0$ such that $\cal N_{j,\bor} = 0$. This implies that $K_jt^\bor.v\neq 0$ for all $v\in V_{aff}^+$ and by lemref{sl2rep}(ii), given $w$ in $W_{tor}$, $w.(K_jt^\bor.v)\neq 0$ whenever $K_jt^\bor.v\neq 0$. Using \corref{sl2rep} assume that $v_0$ is a fixed vector in $V_{aff}^+$ with  $K_jt^\bor.v_0\neq 0$ and weight of $v_0$ is  $\lambda=\lambda|_{\lie h_{aff}} +
\underset{i=2}{\overset{k}{\sum}} m_i\delta_i$, with $r_i-m\leq m_i<r_i$ for $2\leq i\leq k$ and $\lambda|_{\lie h_{aff}}\in P_{aff}^+$.

%We first consider the case when $r_1=0$. 

Let $\btilde$ be the $\Z^{k-1}$-graded Lie subalgebra of $\lie h_\fn\otimes \C[t_1^\pm,\cdots,t_k^\pm]\oplus \cal Z\oplus D_k$ generated by $\lie h_{aff}\otimes \C[t_2,\cdots,t_k]$. %\oplus D_{k-1}$ where $D_{k-1}$ denotes the $\C$-linear span of $d_2,\cdots,d_k$ .
Let  $\cal {V}^\lambda =\cal U(\btilde).v_0$ be the $\btilde$-module generated by $v_0$. Clearly $\btilde$ is a solvable Lie algebra and $\cal V^\lambda$ is a $\btilde$-submodule of $V_{aff}^+$.  Fixing   $$\eta_j = \alpha_{n+1}-\delta_j, \hspace{1cm} \text{and} \hspace{1cm} w_j = r_{\alpha_{n+1}}r_{\eta_j}, \hspace{1cm} \text{for}\ 2\leq j\leq k, $$
and setting $\bow^{\underline\bll} =   w_k^{l_k}\cdots w_3^{l_3}w_2^{l_2}$ for
$\underline\bll=(l_2,\cdots,l_k)\in \Z_+^{k-1}-\{\bold{0}\}$,   we see that 
$$\cal{V}^{\lambda,\bll}:= \cal U(\btilde).\bow^{\underline\bll}(v_0)$$ is a proper submodule of  $\cal V^\lambda$ and the corresponding quotient space
$\cal V^\lambda/\cal V^{\lambda,\bll}$ is a non-zero finite-dimensional module for the solvable Lie algebra $\btilde$. 
Hence by Lie's theorem there exists
$\bar{u}_0\in \cal V^\lambda/\cal V^{\lambda,\bll}$ and a function
$\phi:\btilde \rightarrow \C$ such that 
$$h\otimes t_2^{p_2}\cdots t_k^{p_k}.\bar{u}_0 = \phi(h,\underline\bop)\bar{u}_0,\hspace{1cm} \forall\ h\in \lie h_\aff, \ \underline\bop=(p_2,\cdots,p_k)\in \Z_+^{k-1}.$$ In particular,
$$\begin{array}{l}h\otimes t_j^{p_j}.\bar{u}_0 = \phi(h,(0,\cdots,0,p_j,0\cdots,0))\bar{u}_0,\\
h'\otimes t_2^{p_2}\cdots t_{j-1}^{p_{j-1}}t_{j+1}^{p_{j+1}}\cdots t_k^{p_k}.\bar{u}_0 = \phi(h,(p_2,\cdots,p_{j-1},0,p_{j+1},\cdots,p_k))
\bar{u}_0, \end{array}$$ for $ h, h'\in \lie h_{aff}, p_2,\cdots,p_k\in \Z_+.$ Hence using the fact that the Killing form $(.|.)$ is non-degenerate on $\lie h_\fn$ and
$$[h\otimes t_j^{p_j}, h'\otimes t_2^{p_2}\cdots t_{j-1}^{p_{j-1}}t_{j+1}^{p_{j+1}}\cdots t_k^{p_k}].\bar{u}_0 = (h|h')p_j K_j t^{\underline{\bop}}.\bar{u}_0 $$ we get that $K_j t^{\underline{\bop}}.\bar{u}_0 =0$ for all $\underline{\bop} \in \Z_+^{k-1}$ with $p_j\neq 0$ This implies
that %for all $\underline\bor=(r_2,\cdots,r_k)\in \Z^{k-1}_+$ and $2\leq j\leq k$, the
graded central elements $K_jt^{\underline\bop}$ act trivially on $\bar{u}_0$ or
equivalently there exists a non-zero weight vector $u_0\in \cal V^\lambda$ such that for all $\underline\bop=(p_2,\cdots,p_k)\in \Z^{k-1}_+$ with $p_j\neq 0$ and $2\leq j\leq k$, $$K_jt^{\underline\bop}.u_0\in \cal V^{\lambda,\bll}.$$ By construction, however, for every weight vector $u\in \cal V^{\lambda,\bll}$, 
$$d_j.u \geq  m_j+m l_j,\hspace{1cm} \text{for all}\  j\geq 2.$$ Hence it follows that $K_jt^{\underline\bop}.u_0=0$ for all $\underline\bop\in \Z_+^{k-1}-\{\bold{0}\}$ with $0\leq p_s< m_s+m l_s$ and $2\leq s\leq k$ whenever $p_j \neq 0$. Choosing $\underline\bll\in \Z_+^{k-1}$ appropriately  we thus see that $K_jt^{\underline\bor}.u_0 =0 $ for $2\leq j\leq k$ which is a contradiction to our initial assumption in the case when $\bor=(r_1,\cdots,r_k)\in \Z^{k}$ is such that $r_1=0$. This shows that $\cal N_{j,\bor} \neq 0$ for $2\leq j\leq k$ and $\underline{\bor}\in \Z^{k-1}$ whenever $r_1=0$ and $r_j\neq 0$. 

Consider now the case when $\bor\in \Z^k$ is such that $r_1<0$. Let $V_{aff}^-=\{u\in V: \lie n_{aff}^- \otimes \C[t_2^{\pm1},\cdots,t_k^{\pm1}].u=0\}$. If $V_{aff}^-$ is a non-empty subset of $V$, then by description,  $$h\otimes t^\bor.u=0, \hspace{.35cm} x_\alpha^-\otimes t_i^s.u=0, \hspace{.35cm} x_\alpha^+\otimes t^\bor.u=0$$ for all $u\in V_{aff}^-$ and $h\in \lie h_\fn, \ \alpha\in R_\fn^+$ $2\leq i\leq k$, $s\in \Z$ , whenever $\bor\in \Z^k$ is such that $r_1<0$. Hence $\cal U(\lie h_\fn\otimes \C[t_2^{\pm1},\cdots,t_k^{\pm1}]).u \subset V_{aff}^-$ for all $u\in V_{aff}^-$ and following the same arguments as in step 1 we can conclude that $\cal N_{j,\bor} \neq 0$ for $1\leq j\leq k$ whenever $\bor\in \Z^k$ is such that $r_1<0$. On the contrary suppose that $V_{aff}^- = \emptyset$ and for every $u\in V$ there exists $t^\bor\in \C[t_1^{-1},t_2^{\pm1},\cdots,t_k^{\pm1}]$ such that $h\otimes t^\bor.u \neq 0$ for some $h\in \lie h_\fn$. Then given a non-zero weight vector $v_0\in V_{aff}^+$, there exists $t^\bop\in \C[t_1^{-1},t_2^\pm,\cdots,t_k^\pm]$ such that $-p_1$ is the least positive
integer for which $h\otimes t^\bop.v_0 \neq 0$. Let $$\cal Y^\lambda = \cal U(\btilde^a).v_0 \hspace{1cm}\text{and} \hspace{1cm} \cal Y^{\lambda,\bll} = \cal U(\btilde^a).\bos^{\underline\bll}(h\otimes  t^\bop.v_0),$$ where $\btilde^a$ is the solvable Lie subalgebra of $\lie h_\fn\otimes \C[t_1^\pm,\cdots,t_k^\pm]\oplus \cal Z\oplus D_k$ generated by $\lie h\otimes \C[t_1^{-1},t_2^{a_2},\cdots,t_k^{a_k}]$, where
$a_j = (\sgn p_j) 1$ for $j\geq 2$ and $\bos^{\underline\bll} = s_2^{l_2}\cdots s_k^{l_k}$ is an element of the Weyl group $W_{tor}$ of $\cal T(\lie g)$ for $\underline\bll =(l_2,\cdots,l_k)\in \Z_+^{k-1}$ with $s_j = r_{\alpha_{n+1}}r_{(\alpha_{n+1}-\delta_j)}$ when $a_j>0$ and $s_j = r_{\alpha_{n+1}}r_{(\alpha_{n+1}+\delta_j)}$ when $a_j<0$. Clearly $\cal Y^{\lambda,\bll}$ is a $\btilde^a$-submodule of $\cal Y^\lambda$ and the corresponding quotient $\cal Y^\lambda/\cal Y^{\lambda,\bll}$ is a finite-dimensional module for $\btilde^a$. Thus by Lie's theorem there exists $\bar{\omega}_0\in \cal Y^\lambda/\cal Y^{\lambda,\bll}$ and an algebra homomorphism $\psi:\btilde^a\rightarrow \C$ such that $$h\otimes t^\bor.\bar{\omega}_0 = \psi(h,\bor)\bar{\omega}_0, \hspace{1.5cm}\forall \ h\otimes t^\bor \in \lie h_\fn \otimes \C[t_1^{-1},t_2^{a_2},\cdots,t_k^{a_k}].$$ By the same arguments as above this implies that $K_jt^\bor.\bar{\omega}_0 =0$  for all $t^\bor\in \C[t_1^{-1},t_2^{a_2},\cdots,t_k^{a_k}]$ and $2\leq j\leq k$ whenever $r_j\neq 0$. Choosing $\underline\bll\in\Z_+^{k-1}$ appropriately and repeating the same arguments as above we conclude that $K_jt^\bor$ acts trivially on $V$ for all $t^\bor\in \C[t_1^{-1},t_2^{a_2},\cdots,t_k^{a_k}]$ whenever $r_j\neq 0$. It is now easy to see that modifying the first part of step 2 of the proof appropriately and using the fact that $\sum_{j=1}^k r_j K_jt^\bor =0$, one
can show that $K_jt^\bor$ acts trivially on $V$ for $1\leq j\leq k$ whenever $\bor=(r_1,\cdots,r_k)\in \Z^k$ is such that $r_1<0$ and $r_j\neq 0$.

\noindent {\bf{Step 3}}. To complete the proof of the proposition it remains to show that  $\cal N_{j,\bor}\neq 0$ for $2\leq j\leq k$ with $\bor\in \Z^{k}$ such that $r_j=0$ and $r_1\leq 0$ and consequently $\cal N_{1,\bor} \neq 0$ whenever $\bor\in \Z^k$ is such that $r_1<0$. By \propref{highest.wt}, $V_{aff}^+$ is non-empty and by step 1 of the proof the subspace $\cal U(\lie h\otimes \C[t_2^{\pm1},\cdots,t_k^{\pm1}]\oplus \cal Z).v$ of $V$ is contained in $V_{aff}^+$ for all $v\in V_{aff}^+$. Thus if $\lambda\in P_{aff}^+$ is such that $
h.v = \lambda(h)v$, for $v\in V_{aff}^+$ and $h\in \lie h_\fn,$
then for all $u\in\cal U(\lie h\otimes \C[t_2^{\pm1},\cdots,t_k^{\pm1}]\oplus \cal Z).v$, $u\in V_{aff}^+$ and 
$$h.u = \lambda(h)u, \hspace{.5cm} \forall\  h\in \lie h_\fn.$$ Since $V$ is integrable, using the representation theory of $\lie{sl}_2(\C)$ it follows that for any $\bor\in \Z^k$, and $\alpha\in R_\fn^+$, 
$$(x_\alpha^-\otimes t^{\bor})^{\lambda(\alpha^\vee)+1}.u =0, \hspace{.5cm} \forall\ u\in V_{aff}^+.$$ Given $\bor =(r_1,\cdots,r_k)\in \Z^k$ set   $\bar{\bor}_j=(r_1,\cdots,r_{j-1},0,r_{j+1},\cdots,r_k)$ for $1\leq j\leq k$. By definition $$[x_\alpha^+\otimes t_i^{-r_i}t^{\bar{\bor}_i}, x_\alpha^-\otimes t_i^{r_i}] = -r_i K_i t^{\bar{\bor}_i} +\sum_{2\leq j\leq k, j\neq i,1} r_j K_j t^{\bar{\bor}_i}  .$$ But by step 2 $K_jt^\bor.v=0$ for all $v\in V$ whenever $\bor\in\Z^k$ is such that $r_j\neq 0$. Hence, 
$$[x_\alpha^+\otimes t_i^{-r_i}t^{\bar{\bor}_i}, x_\alpha^-\otimes t_i^{r_i}].u = -r_i K_i t^{\bar{\bor}_i}.u ,\hspace{.75cm} \text{whenever}\ u\in V_{aff}^+.$$ 

For $1\leq i\leq k$, set 
$$D_{\alpha,i}^{\bar{\bor}_i}(u) = \exp (-\sum_{r=1}^\infty \frac{\alpha^\vee \otimes (t^{\bar{\bor}_i})^s+K_i(t^{\bar{\bor}_i})^s}{s} u^s) = \sum_{\ell=1}^\infty D_{\alpha,i}^{\bar{\bor}_i}(\ell)\ u^\ell .$$ 
If $\bor\in \Z^k$ is  such that $r_1=0$ and $r_i\neq 0$ for some $i\geq 2$, then using the integrability condition on a vector $v\in V_{aff}^+\cap V_\lambda$ and the fact that $K_jt^\bor$, $2\leq j\leq k$, acts trivially on $V$ whenever $\bor\in \Z^k$ is such that $r_j\neq 0$ we see that,
$$
\frac{(x_\alpha^+\otimes t_i^{-r_i}t^{\bar{\bor}_i})^{s}}{s!}\frac{(x_\alpha^-\otimes t_i^{r_i})^{s+1}}{(s+1)!}.v = \sum_{\ell=0}^s (x_\alpha^-\otimes t_i^{r_i}(t^{\bar{\bor}_i})^{s-\ell}) D_{\alpha,i}^{\bar{\bor}_i}(\ell).v = 0, \hspace{.35cm} \text{whenever}\ s\geq \lambda(\alpha^\vee).
$$
Applying $x_\alpha^+$ to the above equation we see that for $v\in V_{aff}^+$
$$\alpha^\vee \otimes t_i^{r_i}(t^{\bar{\bor}_i})^{\lambda(\alpha^\vee)}.v =-(\sum_{\ell=1}^{\lambda(\alpha^\vee)} (\alpha^\vee\otimes t_i^{r_i}(t^{\bar{\bor}_i})^{s-\ell}) D_{\alpha,i}^{\bar{\bor}_i}(\ell)).v $$
By Step 2 it follows that  for all $v\in V_{aff}^+\cap V_\lambda$,
$$
[\ h\otimes (t_j^{-r_j})^{(\lambda(\alpha^\vee)-1)}, \alpha^\vee \otimes t_i^{r_i}(t^{\bar{\bor}_i})^{\lambda(\alpha^\vee)}].v= -r_j(h|\alpha^\vee)(\lambda(\alpha^\vee)-1)K_j t_j^{r_j}t_i^{r_i}(\prod_{s\neq 1,i,j}t_s^{r_s})^{\lambda(\alpha^\vee)}.v=0.
$$  Hence
$$ \sum_{\ell=1}^{\lambda(\alpha^\vee)} [\ h\otimes (t_j^{-r_j})^{(\lambda(\alpha^\vee)-1)},(\alpha^\vee\otimes t_i^{r_i}(t^{\bar{\bor}_i})^{(\lambda(\alpha^\vee)-\ell)}) D_{\alpha,i}^{\bar{\bor}_i}(\ell)].v=0, \hspace{.35cm}\forall\ v\in V_{aff}^+.$$
Again using the fact that $K_jt^\bor$ acts trivially on $V$ whenever $\bor\in \Z^k$ is such that $r_j\neq 0$ we deduce from the above equation that for all $v\in V_{aff}^+$,
$$\begin{array}{r}
r_j(\lambda(\alpha^\vee)-1)(h|\alpha^\vee)\big[ K_j(t_i^{r_i}(\prod_{s\neq 1,i,j}t_s^{r_s})^{\lambda(\alpha^\vee)-1})(\alpha^\vee\otimes t^{\bar{\bor}_i} +K_it^{\bar{\bor}_i}) +\\- (\alpha^\vee\otimes t^{\bor}) K_j (\prod_{s\neq 1,i,j} t_s^{r_s})^{\lambda(\alpha^\vee)-1} + \\ 2(\alpha^\vee\otimes t_i^{r_i})(\alpha^\vee\otimes (t^{\bar{\bor}_i})+ K_i(t^{\bar{\bor}_i})) K_j(\prod_{s\neq 1,i,j} t_s^{r_s})^{\lambda(\alpha^\vee)-1}\big].v=0.
 \end{array}$$ Since the bilinear form $(.|.)$ is non-degenerate on $\lie h_\fn$, choosing $h\in \lie h_\fn$ appropriately we see that for any $v\in V_{aff}^+$, 
 $$\begin{array}{r}
 (\alpha^\vee \otimes t^\bor) K_j (\prod_{s\neq 1,i,j} t_s^{r_s})^{\lambda(\alpha^\vee)-1}.v = K_j(t_i^{r_i}(\prod_{s\neq 1,i,j}t_s^{r_s})^{\lambda(\alpha^\vee)-1})(\alpha^\vee\otimes t^{\bar{\bor}_i} +K_it^{\bar{\bor}_i})+\\ 2(\alpha^\vee\otimes t_i^{r_i})(\alpha^\vee\otimes (t^{\bar{\bor}_i})+ K_i(t^{\bar{\bor}_i})) K_j(\prod_{s\neq 1,i,j} t_s^{r_s})^{\lambda(\alpha^\vee)-1}.v
  \end{array}
$$  
 Now applying $h\otimes t_i^{-r_i}$ to the above equation and using first part of the proof and the fact that $K_i$ acts trivially on $V$ for $i\neq 1$ we see that there exists $v_0\in V_{aff}^+$ such that $K_i t^{\bar{\bor}_i}.v_0=0,$ for $i\geq 2$. Since $\bor\in \Z^k$ is arbitrary we see that $K_it^\bor$ acts trivially on V for all $\bor\in \Z^k$ such that $r_1=0=r_i$.
 
Finally we consider the case when $\bor\in \Z^k$ is such that $r_1<0$ and $r_i\neq 0$ for some $2\leq i\leq k.$ Using the integrability condition as above we see that for $v\in V_{aff}^+\cap V_\lambda$,
$$\frac{(x_\alpha^+\otimes t_1^{-{r_1}}t^{\bar{\bor}_1})^s}{s!} \frac{(x_\alpha^-\otimes t_1^{r_1})^{s+1}}{(s+1)!}.v = \sum_{\ell=0}^{s} (x_\alpha^-\otimes t_1^{-r_1}(t^{\bar{\bor}_1})^{s-\ell}) D_{\alpha,1}^{\bar{\bor}_1}(\ell).v =0, \hspace{.35cm} \text{whenever} \ s\geq \lambda(\alpha^\vee).$$ Applying $\ad (h\otimes (t_i^{-r_i})^{(\lambda(\alpha^\vee)-1)})\ad x_{\alpha}^+$ to the above equation and using the first part of the proof we get that for $v\in V_{aff}^+\cap V_\lambda$,
$$\begin{array}{r} 0= K_it_1^{r_1}t_i^{r_i}(\prod_{s\neq 1,i}t_s^{r_s})^{\lambda(\alpha^\vee)}.v = K_it_1^{r_1}(\prod_{s\neq 1,i}t_s^{r_s})^{\lambda(\alpha^\vee)-1}(\alpha^\vee\otimes t^{\bar{\bor}}_1+K_1t^{\bar{\bor}_1}).v \\ - (\alpha^\vee\otimes t^\bor)(K_i(\prod_{s\neq 1,i} t_s^{r_s})^{\lambda(\alpha^\vee)-1}).v\\ +2 (\alpha^\vee\otimes t_1^{r_1})(\alpha^\vee \otimes t^{\bar{\bor}_1}+K_1t^{\bar{\bor}_1})K_i(\prod_{s\neq i,i} t_s^{r_s})^{\lambda(\alpha^\vee)-1}.v,\end{array} $$ which implies that for all $v\in V_{aff}^+\cap V_\lambda$,
$$\begin{array}{r}  (\alpha^\vee\otimes t^\bor)(K_i(\prod_{s\neq 1,i} t_s^{r_s})^{\lambda(\alpha^\vee)-1}).v= K_it_1^{r_1}(\prod_{s\neq 1,i}t_s^{r_s})^{\lambda(\alpha^\vee)-1}(\alpha^\vee\otimes t^{\bar{\bor}_1}+K_1t^{\bar{\bor}_1}).v \\\ +2 (\alpha^\vee\otimes t_1^{r_1})(\alpha^\vee \otimes t^{\bar{\bor}_1}+K_1t^{\bar{\bor}_1})K_i(\prod_{s\neq i,i} t_s^{r_s})^{\lambda(\alpha^\vee)-1}.v.\end{array} $$ Applying $\ad(h\otimes t_i^{-r_i})$ to the above equation and using the fact that for $2\leq i\leq k$, $K_it^\bor$ acts trivially on $V$ whenever $\bor\in \Z^k$ is such that $r_1=0=r_i$ or $r_i\neq 0$ and $r_1 <0$ , we see that $\cal N_{i,\bar{\bor}_i}\cap (V_{aff}^+\cap V_\lambda)\neq 0$. Since for any $\bor\in \Z^k$, $\sum_{j=1}^k r_jK_j t^{\bor} =0 $ it follows that $\cal N_{1,\bor}\neq 0$ whenever $r_1 < 0.$ This completes the proof of the proposition.
\endproof \label{gr.central.action}
It follows from \propref{gr.central.action} that an irreducible $\cal T(\lie g)$-module in $\cal I_\fn^{(m\boe_1)}, m>0$ is in fact a module for the Lie algebra $\lie g_\fn\otimes \C[t_1^{\pm1},\cdots,t_k^{\pm1}] \oplus \cal Z_1 \oplus D_k$, where $\cal Z_1$ is the subspace of $\cal Z$ spanned by the central elements $K_1t^{\bor}$, with $\bor =(r_1,\cdots,r_k)\in \Z^k$  such that $r_1=0$.

\subsection{} \label{A.Lambda} Let $V$ be an irreducible $\cal T(\lie g)$-module in $\cal I_\fn^{(m\boe_1)}$. By \propref{highest.wt}  there exists a non-zero
weight vector $v_0\in V_{aff}^+$ such that $V= \cal U(\cal T(\lie g)).v_0$ and %Given $v_0\in V_{aff}^+$, let $V_0$ be the $\cal T(\lie g)$-submodule of $V$ generated by the vector $v_0$ of weight $\Lambda$. Since $V_0\subseteq V$ is integrable module in $I_{fin}^{(m)}$, $m>0$ $\Lambda\in P_{aff}^+$ and \\
\begin{equation}\label{five.one}
\lie n_{aff}^+\otimes \C[t_2^{\pm1},\cdots,t_k^{\pm1}].v_0=0, \hspace{.35cm} h.v_0 = \Lambda(h).v_0\, \ \forall \ h\in \lie h_{aff},  \hspace{.35cm} (x_i^-)^{\Lambda(\alpha_i^\vee)+1}.v_0 =0,\, \ \text{for}\ i=1,\cdots,n+1.\end{equation} Clearly $V_{aff}^+\subseteq \cal U(\lie h_{aff}\otimes \C[t_2^{\pm1},\cdots,t_k^{\pm1}]).v_0$. For $\alpha\in R_\fn$ and $r\in \Z$ the triple $\{x_\alpha^\pm\otimes t_1^{\pm r}, \alpha^\vee-2r/(\alpha|\alpha)K_1 \}$ is isomorphic to a copy of $\lie {sl}_2(\C)= \C x^+\oplus \C x^-\oplus\C h$. Hence given a monomial $a\in \C[t_2^\pm,\cdots,t_k^\pm]$, the map  $\phi_{\alpha,r}^a:\lie {sl}_2(\C)\otimes \C[t^{\pm 1}]\rightarrow \lie {sl}_2(\C)\otimes \C[a^{\pm 1}]$ given by\\
$$x^\pm\otimes t^s\mapsto (x_{\alpha}^\pm\otimes t^{\pm r})\otimes a^s, \hspace{.35cm} h\otimes t^s\mapsto   (\alpha^\vee-2r/(\alpha|\alpha)K_1)\otimes a^s,\hspace{.45cm} \forall\ s\in \Z,$$\\
defines a Lie algebra homomorphism. For $\beta= \alpha+r\delta_1 \in R^{{re}^{+}}_{aff}$, let $p_{\beta,a}^s = \phi_{\alpha,r}^a(p_\alpha^s)$. Using the homomorphism $\phi_{\alpha,r}^a$ and \lemref{Garland} it is easy to see that for $v_o\in V_{aff}^+$ we have,
\begin{equation}
(x_\beta^+\otimes a)^s(x_\beta^-\otimes 1)^{s+1}.v_0=\underset{l=0}{\overset{s}{\sum}}(x_\beta^-\otimes a^{s-l})p_{\beta,a}^l.v_0 , \label{five.two}
\end{equation}
\begin{equation}\label{five.three}
(x_\beta^+\otimes a)^{s+1}(x_\beta^-\otimes a)^{s+1}.v_0 = p_{\beta,a}^{s+1}.v_0.
\end{equation}
Note that $p_{\beta,a}^l\in \cal U(\lie h_{aff}\otimes\C[t_2^{\pm 1},\cdots,t_k^{\pm 1}])$, hence applying $x_{\beta}^+$ to equation \eqref{five.two} we get $$x_\beta^+(x_\beta^+\otimes a)^s(x_\beta^-\otimes 1)^{s+1}.v_0=\underset{l=0}{\overset{s}{\sum}}(\beta^\vee\otimes a^{s-l})p_{\beta,a}^l.v_0 .$$
Given $\Lambda\in P_{aff}^+$ and $v_0\in V_\Lambda\cap V_{aff}^+$, we thus conclude from \eqref{five.one}, \eqref{five.two} and \eqref{five.three}, that for all $\beta\in R_{aff}^{{re}^+}$ and $\ s\geq {\Lambda(\beta^\vee)+1},$
\begin{equation}\label{five.four}p_{\beta,a}^s.v_0 =0, \hspace{.65cm}\ \beta^\vee \otimes a^s.v_0 +( \underset{l=1}{\overset{s}{\sum}}(\beta^\vee\otimes a^{s-l})p_{\beta,a}^l).v_0=0.\end{equation}

\vspace{.35cm}

For a fixed $\Lambda\in P_{aff}^+$, let $I_\Lambda$ be the ideal of $\cal U(\lie h_{aff}\otimes \C[t_2^{\pm1},\cdots,t_k^{\pm1}])$ generated by the elements:
$$h-\Lambda(h), \hspace{.25cm} \forall \ h\in \lie h_{aff}, \hspace{.5cm} K_jt^\bom \hspace{.25cm} 2\leq j\leq k, \ \text{and all monomials} \ t^\bom \ \text{in} \ \C[t_2^{\pm1},\cdots,t_k^{\pm1}]$$  
$$p_{\alpha_i,a}^s, \hspace{.35cm} \text{for}\ |s|>\Lambda(\alpha_i^\vee), \ i=1,\cdots, n+1, \ \text{and all monomials} \ a\ \text{in} \ \C[t_2^{\pm1},\cdots,t_k^{\pm1}],$$
$$\sum_{l=0}^s (\alpha_i^\vee \otimes a^{s-l})p_{\alpha_i,a}^l, \hspace{.35cm} \text{for}\ |s|>\Lambda(\alpha_i^\vee), \ i=1,\cdots, n+1, \ \text{and all monomials} \ a\ \text{in} \ \C[t_2^{\pm1},\cdots,t_k^{\pm1}].$$
Let $$\bold{A}_\Lambda =  \cal U(\lie h_{aff}\otimes \C[t_2^{\pm1},\cdots,t_k^{\pm1}])/I_\Lambda.$$
%Since $\C[t_2^\pm,\cdots,t_k^\pm]$ is finitely generated, and by 
Using \eqref{five.four} it is easy to see that the algebra $\bold{A}_\Lambda$ is generated by the elements of the set $$\{\alpha_i^\vee \otimes t_2^{s_2}t_3^{s_3}\cdots t_k^{s_k} : |s_l|\leq \Lambda(\alpha_i^\vee), \ \text{for} \ 2\leq l\leq k,\ i=,2,\cdots, n+1\}.$$ Hence $\bold{A}_\Lambda$ is a finitely generated commutative algebra. Further it  follows from  \eqref{five.four} that given a weight vector $v_0\in V_{aff}^+\cap V_{\Lambda},$ the $\lie h_{aff}\otimes \C[t_2^{\pm1},\cdots,t_k^{\pm1}]$-submodule of $V_{aff}^+$ generated by $v_0$ is a quotient of the algebra $\bold{A}_\Lambda.$ Thus the  irreducible $\lie h_{aff}\otimes \C[t_2^{\pm1},\cdots,t_k^{\pm1}]$-submodules of $V_{aff}^+$ generated by a non-zero vector of weight $\Lambda\in P_{aff}^+$ are in one-to-one correspondence with the maximal ideals of $\bold{A}_\Lambda$, that is, the irreducible $\lie h_{aff}\otimes \C[t_2^{\pm1},\cdots,t_k^{\pm1}]$-submodules of $V_{aff}^+$ generated by a non-zero vector of weight $\Lambda\in P_{aff}^+$ are in one-to-one correspondence with the set of graded algebra homomorphisms from $\bold{A}_\Lambda$ to $\C[t_2^{\pm1},\cdots,t_k^{\pm1}]$.

\subsection{} Set $\cal L^{c}(\lie g) : =\lie g_\fn\otimes \C[t_1^{\pm1},\cdots,t_k^{\pm1}] \oplus \cal Z_1 \oplus \C d_1$. Clearly $\lie h_fin \oplus \C K_1 \oplus \C d_1 = \lie h_{aff}$ is the Cartan subalgebra of $\cal L^c(\lie g).$ Let $\ev({\bf{1}}):\C[t_2^{\pm1},\cdots,t_k^{\pm1}]\rightarrow \C$ be the evaluation map defined by $\ev({\bf{1}})(t^\bom) =1$ for $t^\bom\in \C[t_2^{\pm1},\cdots,t_k^{\pm1}]$.

Let $V$ be an irreducible $\cal T(\lie g)$-module in $\cal I_\fn^{(m\boe_1)}$ and $v_0\in V_{aff}^+$ be a non-zero vector of weight $\Lambda.$ Then following our discussion in (\ref{A.Lambda}) there exists a graded algebra homomorphism $\phi:\bold{A}_\Lambda \rightarrow \C[t_2^{\pm1},\cdots,t_k^{\pm1}]$ such that
$ Y.v_0 = \phi(Y).v_0,$ for all $Y\in \bold{A}_\Lambda.$
Given $\phi:\bold{A}_\Lambda\rightarrow \C[t_2^{\pm1},\cdots,t_k^{\pm1}]$, let $\bar\phi :\bold{A}_\Lambda \rightarrow \C$ be the algebra homomorphism defined by  the composition of functions $\bar\phi := \ev({\bf{1}})\circ\phi$ and let
$W_{\bar{\phi}} := \cal U(\cal L^c(\lie g)).v_0$ be the integrable $\cal L^c(\lie g)$ module generated by the vector $v_0$ such that
$Y.v_0 = \ev({\bf{1}}).\phi(Y)v_0 = \bar{\phi}(Y).v_0,$ for all $Y\in \bold{A}_\Lambda.$  Since $$\lie n_{aff}^+\otimes \C[t_2^{\pm1},\cdots,t_k^{\pm1}].v_0 =0$$ and by definition the highest weight space $\cal U(\lie h_\fn \otimes \C[t_2^{\pm1},\cdots, t_k^{\pm1}]\oplus \cal Z_1).v_0$  of $W_{\bar{\phi}}$ is one-dimensional, $W_{\bar{\phi}}$ has a unique irreducible quotient $V_{\bar{\phi}}$.

\begin{lem} Let $\Lambda\in P_{aff}^+$. Given an algebra homomorphism $\psi: \bold{A}_\Lambda \rightarrow \C $, let $W_\psi$ be the integrable $\cal L^c(\lie g)$-module generated by a vector $v$ of weight $\Lambda$ such that $$\lie n_{aff}^+\otimes \C[t_2^{\pm1},\cdots,t_k^{\pm1}].v =0, \hspace{.5cm} Y. v=\psi(Y)v\ \hspace{.25cm} \forall \ Y\in \bold{A}_\Lambda, \hspace{.5cm} K_1.v = mv.$$ Then $W_\psi$ is an integrable $\cal L^c(\lie g)$-module with finite-dimensional weight spaces. 
  \end{lem}
\proof $W_\psi$ is an integrable $\cal L^c(\lie g)$-module and hence can be written as $$W_{\psi} =\underset{\mu\in \lie h_{aff}^*}{\oplus} W_{\psi}^\mu, \hspace{.75cm} \text{where}\ W_\psi^\mu =\{u\in W_\psi: h.u = \mu(h)u, \ \forall\ h\in \lie h_{aff}\}.$$ Let $P(W_\psi):= \{\mu\in \lie h_{aff}^*: W_\psi^\mu\neq 0\}$. Since $W_\psi$ is an integrable $\cal L^c(\lie g)$-module with highest weight $\Lambda\in P_{aff}^+$, an element $\mu\in P(W_\psi)$ is of one of the following forms:
\item[(i).]$\mu = \Lambda -r\delta_1$, $r\in \Z_{\geq 0}$. 
\item[(ii).] $\mu= \Lambda +\eta-r\delta_1$, $\eta\in Q$ and $r\in \Z_{\geq 0}$.
  
\noindent We prove that in each of the cases $W_\psi^\mu$ is spanned by a finite set of elements and hence is finite-dimensional.

\noindent {\textbf{Case (i)}}. Let $\mu = \Lambda-r\delta_1$ with $r\in \Z_{\geq 0}$. By definition $\dim W_\psi^\Lambda =1$. Suppose now that $r>0$. Then any
element in $W_\psi^{\Lambda-r\delta_1}$ is of the form $(h_1\otimes t_1^{-r_1}t^{\bar{\bor}_1})(h_2\otimes t_1^{-r_2}t^{\bar{\bor}_2})\cdots(h_s\otimes t_1^{-r_s}t^{\bar{\bor}_s}).v$ with $r_i\in \Z_{>0}$ such that $\sum_{i=1}^s r_i =r$, $t^{\bar{\bor}_i}\in \C[t_2^{\pm1},\cdots,t_k^{\pm1}]$ and $h_i\in \lie h_\fn$ for $1\leq i\leq s$. Since
$\Delta_\fn^\vee:=\{\alpha_i^\vee: 1\leq i\leq n\}$ is a basis of $\lie h_\fn$ it
follows that any element in $W_\psi^{\Lambda-r\delta_1}$ is spanned by elements of
the form
\begin{equation}\label{h.r}
(\alpha_{i_1}^\vee\otimes t_1^{-r_1}t^{\bar{\bor}_1})(\alpha_{i_2}\otimes t_1^{-r_2}t^{\bar{\bor}_2})\cdots(\alpha_{i_s}\otimes t_1^{r_s}t^{\bar{\bor}_s}).v,
\end{equation}
with $r_i\in \Z_{>0}$ such that $\sum_{j=1}^s r_j =r$, $t^{\bar{\bor}_j}\in \C[t_2^{\pm1},\cdots,t_k^{\pm1}]$ and $\alpha_{i_j} \in \Delta_\fn^\vee$ for $1\leq j\leq s.$ Observe now that by putting $a=t_1^{-r_1}t_2^{r_2}\cdots t_{i-1}^{r_{i-1}}t_{i+1}^{r_{i+1}}\cdots t_k^{r_k}$  and using Garland's equation (\lemref{Garland}) we get 
 \begin{equation}(x_{\alpha_j}^+\otimes a^{-1}t_i)^{\Lambda(\alpha_j^\vee)} (x_{\alpha_j}^-\otimes a)^{\Lambda(\alpha_j^\vee)+1}.v = \sum_{j=0}^{\Lambda(\alpha_j^\vee)} x_{\alpha_j}^- \otimes a (t_i)^{\Lambda(\alpha_j^\vee)-j} p^j_{\alpha_j,t_i}.v =0.
 \label{beta.1.a}\end{equation}
Applying $x_{\alpha_j}^+\otimes t_i^p$ to \eqref{beta.1.a} we get
\begin{equation} \label{h.r.bor}
  - \alpha_j^\vee \otimes a (t_i)^{\Lambda(\alpha_j^\vee)+p}=  \sum_{\ell=1}^{\Lambda(\alpha_j^\vee)} \alpha_j^\vee \otimes a (t_i)^{\Lambda(\alpha_j^\vee)-\ell+p} p^\ell_{\alpha_j,t_i}.v , \hspace{.5cm} \forall \ p\in \Z.\end{equation} As $K_it^\bor.v =0 $ for all $1\leq i\leq k$ whenever $\bor\in \Z^k$ is such that $r_1\neq 0$ it follows from \eqref{h.r.bor} and the fact that $W_\psi^{\Lambda-r\delta_1}$ is spanned by elements of the form \eqref{h.r}, that $W_\psi^{\Lambda-r\delta_1}, r\in \Z_{>0}$ lies in that span of elements in the set
$$\begin{array}{r}S_{\Lambda-r\delta_1}=\{(\alpha_{i_1}^\vee\otimes t_1^{-r_1}t^{\bar{\bor}_1})(\alpha_{i_2}\otimes t_1^{-r_2}t^{\bar{\bor}_2})\cdots(\alpha_{i_s}\otimes t_1^{r_s}t^{\bar{\bor}_s}).v: \alpha_{i_j}\in \Delta_\fn^\vee,  r_j\in \Z_{>0}, r=\sum_{j=1}^sr_j,\\ \bar{\bor}_j =(0,r_2^j,\cdots,r_k^j)\ \text{with} \ |r^j_t|<\Lambda(\alpha_{i_j})\ 1\leq j\leq s,\ 2\leq t\leq k\}\end{array}.$$
Since $\Delta_\fn^\vee$ is finite and the number of partitions of a given positive number is finite it follows the length of any element in $W_\psi^{\Lambda-r\delta_1}$ is less than $r$. Hence the set $S_{\Lambda-r\delta_1}$ is finite implying that $W_\psi^{\Lambda-r\delta_1}$ is finite-dimensional.

\noindent {\textbf{Case (ii.)}}. Let $\mu = \Lambda +\eta-r\delta_1$ with $r\in \Z_{> 0}$ and $\eta\in Q$. Note that the elements in $W_\psi^{\Lambda+\eta-r\delta_1}$ are of the form
\begin{equation}\label{eta.r}
  (Y_1\otimes t_1^{-r_1}t^{\bar{\bor}_1})(Y_2\otimes t_2^{-r_2}t^{\bar{\bor}_2})\cdots(Y_s\otimes t_1^{-r_s}t^{\bar{\bor}_s}).v,\end{equation}
with $Y_j\in \lie g_\fn$ for $1\leq j\leq s$, $r_j <0$ whenever $Y_j\in \lie n_\fn^+\oplus \lie h_\fn$ and $r_j\leq 0$ for $Y_j\in \lie n_\fn^-.$ %Using the Lie bracket operation in $\cal L^c(\lie g)$ it is easy to see that any element of the form \eqref{eta.r} can be written as a linear combination of elements of the form \begin{equation}\label{eta.r}
 % (Y_1\otimes t_1^{-r_1}t^{\bar{\bor}_1})(Y_2\otimes t_2^{-r_2}t^{\bar{\bor}_2})\cdots(Y_s\otimes t_1^{-r_s}t^{\bar{\bor}_s}).v,\end{equation}
We prove by induction on $s$ that any such element is in the span of the elements 
$$(Y_1\otimes t^{\bom_1})\cdots(Y_\ell\otimes t^{\bom_\ell}).v $$ where for all $1\leq t\leq \ell$, $ Y_t\in \lie g_\fn$ and $\bom_t=(m^t_1,m^t_2,\cdots,m_k^t) \in\Z^k$ is such that $m_1^t\leq 0$ and for $2\leq j\leq k$, $|m_j^t| < \underset{1\leq i\leq n}{\max}\Lambda(\alpha_i^\vee)$  if $Y_t\in \lie h_\fn$ and $|m_j^t| < \Lambda(\gamma_{r_t}^\vee)$  if $Y_t \otimes t_1^{-r_t}$ is a real root vector of $\lie g_{aff}$ of weight $\gamma -r_t\delta_1$. 
%Applying induction on $s$ and using arguments on the partition function on $\N$
%we will show that $W_\psi^{\Lambda+\eta-r\delta_1}$ is finite-dimensional for $\eta\in Q$ and $r<0$.

Suppose $s=1$. Then $\eta\in R_\fn$, $r_1=r$ and  $Y_1\otimes t_1^{-r} \in \lie n_{aff}^-$ is of the form $x_\eta \otimes t_1^{-r}$. Hence there exists $y_\eta\otimes t_1^{r}\in \lie n_{aff}^+$ such that the triple $\lie{sl}_2(\eta,-r)=\{y_\eta\otimes t_1^{r},x_\eta\otimes t_1^{-r}, \eta_{r}^{\vee}:=[y_\eta\otimes t_1^{r},x_\eta\otimes t_1^{-r}]\}$ is isomorphic to $\lie{sl}_2(\C)$. Setting
$a=t_2^{r_2}\cdots t_{i-1}^{r_{i-1}}t_{i+1}^{r_{i+1}}\cdots t_k^{r_k}$ and using Garland's equation we get
$$(y_\eta\otimes t_1^r a^{-1}t_i)^{\Lambda(\eta^\vee_r)}(x_{\eta}\otimes t_1^{-r}a)^{\Lambda(\eta_r^\vee)+1}.v = \sum_{j=0}^{\Lambda(\eta_r^\vee)} x_\eta\otimes t_1^{-r}at_i^{\Lambda(\eta_r^\vee)-j}p_{\eta+r\delta_1,t_i}^j.v$$ Applying $\eta^\vee\otimes t_i^p$ to the above equation we see that for any real root $\gamma \in R_{aff}^+$, the element $x_\gamma^-\otimes t^{\bar{\bor}}.v$, with $\bar{\bor} = (0,r_2,\cdots,r_k)$, lies in the span of $x_\gamma^-\otimes t_2^{s_2}\cdots t_k^{s_k}$ where $-\Lambda(\gamma^\vee)<s_j<\Lambda(\gamma^\vee)$, $2\leq j\leq k$.
If an element in $W_\psi^{\Lambda+\eta-r\delta_1}$ is of the form $(Y_1\otimes t_1^{-r_1}t^{\bar{\bor}_1})(Y_2\otimes t_1^{r_2}t^{\bar{\bor}_2}).v$ with $Y_1\in \lie n_\fn^+\oplus \lie n_\fn^-$ and $Y_2\in \lie h_\fn$, then using the Lie bracket operation in $\cal L^c(\lie g)$, the relations \eqref{h.r.bor} and induction step $s=1$ it follows such elements lie in the linear span of elements of the form
\begin{equation}(Y_1\otimes t_1^{-r_1}t^{\bar{\bom}_1})(Y_2\otimes t_1^{-r_2}t^{\bar{\bom}_2}).v, \label{Y.1.2}\end{equation}
where for $2\leq j\leq k$, $|m^2_j|< \underset{1\leq i\leq n}{\max}\Lambda(\alpha_i^\vee)$ and $|m^1_j|< \Lambda(\eta_{r_1}^\vee)$. Since the number of pairs $(r_1,r_2)\in \Z_{\geq 0}$ such that $r_1+r_2 =r$ is finite, it follows that every element in $W_\psi^{\Lambda+\eta-r\delta_1}$ of the form \eqref{Y.1.2} lies in the span of a finite set of elements.  Suppose now that the result is true when $1\leq s<t $ and consider an element in $W_\psi^{\Lambda+\eta-r\delta_1}$ of length $t$. Clearly
$$\begin{array}{r}(Y_1\otimes t_1^{-r_1}t^{\bar{\bor}_1})(Y_2\otimes t_1^{-r_2}t^{\bar{\bor}_2})\cdots(Y_s\otimes t_1^{-r_s}t^{\bar{\bor}_s}).v = (Y_2\otimes t_1^{-r_2}t^{\bar{\bor}_2})(Y_1\otimes t_1^{-r_1}t^{\bar{\bor}_1})\cdots(Y_s\otimes t_1^{-r_s}t^{\bar{\bor}_s}).v \\+ [Y_1\otimes t_1^{-r_1}t^{\bar{\bor}_1}, Y_2\otimes t_1^{-r_2}t^{\bar{\bor}_2}](Y_3\otimes t_1^{-r_3}t^{\bar{\bor}_3})\cdots(Y_s\otimes t_1^{-r_s}t^{\bar{\bor}_s}).v.\end{array}$$ Since $ [Y_1\otimes t_1^{-r_1}t^{\bar{\bor}_1}, Y_2\otimes t_1^{-r_2}t^{\bar{\bor}_2}] \in \lie g_\fn \otimes t_1^{-r_1-r_2}t^{\bar{\bor}_1+\bar{\bor}_2}$, applying an induction argument on the length $t$ identical to the one used in \cite[Proposition 1.2(ii)]{CPweyl} we see that elements of the form \eqref{eta.r} are spanned by elements of the desired form. Since $R_\fn$ is finite, for a fixed $\eta\in Q$ there can exist only finitely many tuples $(\beta_1,\beta_2,\cdots,\beta_t)\in (R_\fn)^t$ such that $\sum_{j=1}^t \beta_j =\eta$, further the number of partitions of a given positive number is finite. Hence given $\eta\in Q$ and $r\in \N$ there exists a positive integer $N_{\eta,r}$ such that every element in $W_{\psi}^{\Lambda+\eta-r\delta_1}$ is of length less than equal to $N_{\eta,r}$. This shows that the set
$$\begin{array}{r}S_{\Lambda+\eta-r\delta_1}=\{ (Y_1\otimes t^{\bom_1})\cdots(Y_\ell\otimes t^{\bom_\ell}).v: Y_j\in \lie g_\fn \ \text{with}\ \underset{j=1}{\overset{\ell}{\sum}}\wt_{\fn}(Y_j)=\eta,  r_j\geq 0 \ \text{with}\ \underset{j=1}{\overset{\ell}{\sum}} r_j =r,\\ \bar{\bor}_j=(0,r^j_2,\cdots,r^j_k)\ \text{with}\ |r_p^j|<\underset{1\leq i\leq n}{\max} \Lambda(\alpha_i^\vee)\ \text{if}\ Y_j\in \lie h_\fn,\\ |r_p^j|\leq \Lambda(\gamma_{r_j}^\vee)\ \text{if}\ \wt_{aff}(Y_j\otimes t_1^{-r_j}) = \gamma-r_j\delta_1\ \forall \ 2\leq p\leq k  \}.\end{array}$$ is finite and since every element in $W_\psi^{\Lambda+\eta-r\delta_1}$ lies in the span of elements from $S_{\Lambda+\eta-r\delta_1}$ it follows that $W_\psi^{\Lambda+\eta-r\delta_1}$ is finite dimensional for $\eta\in Q$ and $r\geq 0$. 
\endproof

We thus conclude that  given an irreducible $\cal T(\lie g)$-module $V$ one can associate with it an unique irreducible $\cal L^c(\lie g)$-module having finite-dimensional weight spaces. This leads us to the study of the
irreducible integrable  representations of $\cal L^c(\lie g)$ having
finite-dimensional weight spaces.

\subsection{} The following result leads towards the classification of the irreducible integrable $\cal L^c(\lie g)$-modules having finite-dimensional weight spaces.  \label{Ch.Rmn.Thm}

\begin{prop} Let $\cal V$ be an integrable irreducible representation of\,  $\cal L^c(\lie g)$ with finite-dimensional weight spaces. Suppose $\cal V$ is generated by a vector $\bold{v}$ such that 
$$\lie n_{aff}^+\otimes \C[t_2^{\pm1},\cdots,t_k^{\pm1}].\bold{v}=0, \hspace{.5cm} h.\bold{v} = \Lambda(h)\bold{v}, \hspace{.25cm}\forall\ h\in \lie h_{aff}, \hspace{.5cm} K_1.\bold{v}=m\bold{v},$$ for $\Lambda\in P_{aff}^+$. Then %$L_{k-1}(\lie g_{aff}')$ 
  $\cal V$ is isomorphic to an irreducible module for a
  finite direct sum of affine Kac Moody Lie algebras $\lie g_{aff}$.
%there exists finitely many maximal ideals  maximal ideals $\cal M_1,\cdots, \cal M_r$ of $\C[t_2^\pm,\cdots,t_k^\pm]$ such that $$\ker \psi = (\lie g\otimes \C[t_1^\pm]\oplus \C K_1)\otimes (\cap_{i=1}^r \cal M_r).$$
\end{prop}
\proof Let $\psi: \cal L^c(\lie g) \rightarrow \End \cal{V}$ be an irreducible integrable representation of of $\cal L^c(\lie g)$ having finite-dimensional weight spaces such that $\psi(K_1) = m\id_{\bar{V}}$ for $m\in \Z_{>0}$. Clearly such a module is $d_1$-graded, therefore $\ker \psi$ is an ideal of $\lie g_\fn\otimes \C[t_2^{\pm1},\cdots,t_k^{\pm1}] \oplus \cal Z_1$.  By definition $\cal Z_1$ is in the center of $\lie g_\fn\otimes \C[t_1^{\pm1},\cdots,t_k^{\pm1}] \oplus \cal Z_1$ . Hence the adjoint representation of $\lie g_\fn\otimes \C[t_2^{\pm1},\cdots,t_k^{\pm1}] \oplus \cal Z_1$ can be identified with a representation of $\lie g_\fn\otimes \C[t_2^{\pm1},\cdots,t_k^{\pm1}] $. This implies that any ideal of $\lie g_\fn\otimes \C[t_2^{\pm1},\cdots,t_k^{\pm1}] \oplus \cal Z_1$  (which is a module of $\lie g_\fn\otimes \C[t_2^{\pm1},\cdots,t_k^{\pm1}] \oplus \cal Z_1$ with respect to the adjoint representation) is an ideal of the Lie algebra $\lie g_\fn\otimes \C[t_2^{\pm1},\cdots,t_k^{\pm1}]$. In particular $\ker \psi$ is an ideal of $\lie g_\fn\otimes \C[t_2^{\pm1},\cdots,t_k^{\pm1}]$ and by \lemref{roots} there exists an ideal $S$ of $\C[t_2^{\pm1},\cdots,t_k^{\pm1}]$ such that $\ker \psi = \lie g_\fn \otimes S.$ If $\cal Z_1^S$ is the central extension of the Lie algebra $\lie g_\fn\otimes (\C[t_1^{\pm1},t_2^{\pm1},\cdots,t_k^{\pm1}]/S)$, then it follows that the irreducible $\cal L^c(\lie g)$-module $\cal V$ is a faithful representation of the Lie algebra $\lie g_\fn\otimes \C[t_1^{\pm1},t_2^{\pm1},\cdots,t_k^{\pm1}] \oplus \cal Z_1^S\oplus \C d_1.$ Hence  the action of the subalgebra $\lie h_{aff} \otimes (\C[t_2^{\pm1},\cdots,t_k^{\pm1}])$ of
$\cal U(\cal L^c(\lie g))$ on $\bold{v}$ is given by the restriction of the map $\psi$ to the Lie algebra $\lie h_{aff} \otimes (\C[t_2^{\pm1},\cdots,t_k^{\pm1}]/S)$.

On the other hand since $\cal V$ is irreducible, it follows from \ref{A.Lambda}
that the $\lie h_{aff}\otimes \C[t_2^{\pm1},\cdots,t_k^{\pm1}]$-module generated by
$\bold{v}$ corresponds to a maximal ideal of the commutative algebra
$\bold{A}_\Lambda,$ or equivalently to an algebra homomorphism
$\phi :  \bold{A}_\Lambda\rightarrow \C.$ Thus for any $h\in \lie h_{aff}$ and
$a\in \C[t_2^{\pm1},\cdots,t_k^{\pm1}]$ we have, 
$$h\otimes a.\bold{v} = \phi(h\otimes a).\bold{v} = \psi(h\otimes a).\bold{v} =
h\otimes \bar{a}.v_0, $$ where $\bar{a}$ is the image of $a$ in
$\C[t_2^{\pm1},\cdots,t_k^{\pm1}]/S$. Since $\phi(h\otimes a)\in \C$, we conclude
that the ideal $S$ is equal to the intersection of distinct maximal ideals
$\{\cal M_j\}_{j\in J}$ of $\C[t_2^{\pm1},\cdots,t_k^{\pm1}]$. By Chinese
remainder theorem it thus follows that an irreducible integrable
$\cal L^c(\lie g)$-module with finite-dimensional weight spaces  is in fact an
irreducible representation of a Lie algebra of the form
$$
\big(\oplus_{j\in J}(\lie g_\fn\otimes \C[t_1^{\pm1}]\oplus \C K_1) \otimes
(\C[t_2^{\pm1},\cdots,t_k^{\pm1}]/\cal M_j) \big) \oplus \C d_1$$
whenever it is generated by a highest weight vector of weight $\Lambda\in P_{aff}^+$ with $\Lambda(K_1)\neq 0$.

 It is well known that an irreducible representation of 
$\underset{j\in J}{\oplus}(\lie g_\fn\otimes \C[t_1^{\pm1}]\oplus \C K_1)\otimes (\C[t_2^{\pm1},\cdots,t_k^{\pm1}]/\cal M_j)\oplus C d_1$
 is a tensor product of irreducible representations of the Lie algebras 
 $$\{(\lie g_\fn\otimes \C[t_1^{\pm1}]\oplus \C K_1)\otimes
 (\C[t_2^{\pm1},\cdots,t_k^{\pm1}]/\cal M_j) \oplus \C d_1: j\in J\}.$$ Hence for
 each $j\in J$ there exists an irreducible representation $V_j$ of
 $(\lie g_\fn\otimes \C[t_1^{\pm1}]\oplus \C K_1)\otimes (\C[t_2^{\pm1},\cdots,t_k^{\pm1}]/\cal M_j) \oplus \C d_1$  such that $\cal V = \otimes_{j\in J} V_j$.  As a
 consequence the highest weight vector $v_0\in V$ is of the form
 $$\bold{v}= \underset{j\in J}{\otimes} v_j, $$ where  $v_j$ is a weight vector of
 $V_j$ for $j\in J$. Since the maximal ideals are distinct, it is easy
 to see that 
 $$\lie n_{aff}^+ \otimes \C[t_2^{\pm1},\cdots,t_k^{\pm1}].\bold{v} =0, \hspace{.75cm} \text{if and
 only if}\hspace{.75cm} \lie n_{aff}^+ \otimes \C[t_2^{\pm1},\cdots,t_k^{\pm1}].v_j =0, \hspace{.35cm} \
 \forall\ j\in J.$$
 Hence for each $j\in J$, the vector $v_j\in V_j$ is of
 weight $\mu_j$ where $\mu_j\in P_{aff}^+$ is such that $\mu_j(K_1)>0$ and for
 all $h\in \lie h_{aff}$,
 $$h.\bold{v} =\Lambda(h)\bold(v) = (\sum_{j\in J}\mu_j(h)) \otimes_{j\in J} v_j.$$
 In particular, $\sum_{j\in J} \mu_j(\alpha_{n+1}^\vee) = m <\infty. $ Hence the
 set $J$ is finite which completes the proof of the proposition. 
\endproof
\label{non.graded}

It follows directly from \propref{non.graded} that given an algebra
homomorphism $\phi:\bold{A}_\Lambda \rightarrow \C$ one can uniquely associate with it a finitely supported function from $\max \C[t_2^{\pm1},\cdots,t_k^{\pm1}]\rightarrow P_{aff}^+$ that maps a maximal ideal $\cal M_j$ of $\C[t_2^{\pm1},\cdots, t_k^{\pm1}]$ to the weight of the highest weight vector of the irreducible
$(\lie g_\fn\otimes \C[t_1^{\pm1}]\oplus \C K_1) \otimes
(\C[t_2^{\pm1},\cdots,t_k^{\pm1}]/\cal M_j \oplus \C d_1$ module under $\phi$.   
%We have seen that By \propref{highest.wt} any irreducible module $V$ in $\cal I_\fn^{(m\boe_1)}$ is generated by a vector $v_0\in V$ such that $\lie n_{aff}^+\otimes\C[t_2^{\pm1},\cdots,t_k^{\pm1}].v_0=0$. Using \propref{aff.rep} it thus follows that such a vector $v_0$ is of the form $\otimes_i v_{_{\Lambda_i}}$, where for each $i$, $\Lambda_i\in P_{aff}^+$ and $v_{_{\Lambda_i}}$ is the highest weight vector that generates the standard module $X(\Lambda_i)$. With this in mind we now proceed to determine the isomorphism classes of the irreducible objects in $\cal I_\fn^{(m\boe_1)}$. 

\vspace{.25cm}

\subsection{} Given $M\in \max \C[t_2^{\pm1},\cdots,t_k^{\pm1}]$ let $\ev_M: \C[t_2^{\pm1},\cdots,t_k^{\pm1}]\rightarrow \C$ be the evaluation map at the point in $(\C^\ast)^{k-1}$ corresponding to the maximal ideal $M$ of $\C[t_2^{\pm1},\cdots,t_k^{\pm1}]$.

\label{six.five}

Let $\Pi$ be the monoid of finitely supported functions $\pi:\max \C[t_2^{\pm1},\cdots,t_k^{\pm1}] \rightarrow P_{aff}^+$. For $\pi,\pi'\in \Pi$ and $M\in \max \C[t_2^{\pm1},\cdots,t_k^{\pm1}]$  define
$$(\pi+\pi')(M)= \pi(M)+\pi'(M),\hspace{1cm} \supp(\pi)=\{M\in \max \C[t_2^{\pm1},\cdots,t_k^{\pm1}]: \pi(M)\neq 0\},$$ $$ \wt(\pi)=\sum_{M\in\supp(\pi)} \pi(M).$$
%and for $\Lambda\in P_{aff}^+$, set $$\Pi_\Lambda = \{\pi\in \Pi: \wt(\pi)=\Lambda\}.$$ 
For $\pi\in \Pi$, let  $M_1,M_2,\cdots,M_l$ is an enumeration of $\supp(\pi)$ and let 
$X_\pi = {\otimes_{i=1}^l} X(\pi(M_i)), $
be the $\lie g_\fn\otimes \C[t_1^{\pm1},\cdots,t_k^{\pm1}]\oplus \cal Z\oplus\C d_1$-module in which the action of the Lie algebra is defined as follows:
$$Y\otimes f.  v_1\otimes \cdots \otimes v_l = \sum_{i=1}^l ev_{M_i}(f) v_1\otimes\cdots\otimes Y.v_i\otimes\cdots\otimes v_l,$$ where $Y\in \lie g_\fn\otimes \C[t_1^{\pm1}]\oplus \C K_1 \oplus \C d_1$ and $f\in \C[t_2^{\pm1},\cdots,t_k^{\pm1}]$. 
Let $$L(X_\pi) := X_\pi \otimes \C[t_2^{\pm1},\cdots,t_k^{\pm1}], $$
 be the $\cal T(\lie g)$-module on which the action of $\cal T(\lie g)$ is given as follows:$$Y\otimes f. (w\otimes f') = (Y\otimes f.w)\otimes ff', \hspace{1.5cm} K_jt^\bom. (w\otimes f') =0, \hspace{.35cm}  \forall\  2\leq j\leq k, \bom\in \Z^k,$$
$$d_i. (w\otimes f') = w\otimes d_i(f'), \hspace{.35cm} \text{for}\ 2\leq i\leq k, \hspace{1.75cm} d_1. w\otimes f' = d_1(w)\otimes f',$$ 
 for $Y\in \lie g_\fn\otimes \C[t_1^{\pm1}]\oplus \C K_1 \oplus \C d_1,$ $f,f'\in \C[t_2^{\pm1},\cdots,t_k^{\pm1}]$ and $w\in X_\pi$. For $M\in \supp(\pi)$, let $v_M$ be the highest weight vector of $X(\pi(M))$  and let $v_\pi:= \otimes_{i=1}^l v_{M_i}.$ 
Suppose $\wt(\pi)= \Lambda$. Then with respect to the action of $\cal T(\lie g)$ on $L(X_\pi)$,  it is clear that the $\bold{A}_{\Lambda}$-module generated by $v_\pi$ is $\Z^{k-1}$-graded. Hence setting
$$A_\Lambda^\pi := \bold{A}_\Lambda.v_\pi,$$ we can write $$A_{\Lambda}^\pi = \underset{\bom\in \Z^{k-1}}{\oplus} A_\Lambda^\pi[\bom].$$
Let $$G_\pi:=\{\bom =(m_2,\cdots,m_k)\in \Z^{k-1}: A_\Lambda^\pi[\bom] \neq 0\}.$$  

\begin{lem} Let $\Lambda\in P_{aff}^+$ and $\pi\in \Pi$ be such that $\wt(\pi)=\Lambda$. Then the set $G_\pi$ is a subgroup of $\Z^{k-1}$ of rank $k-1$. 
\end{lem}
\proof Given $\pi\in \Pi$ with $\wt(\pi)=\Lambda\in P_{aff}^+$, there exists some  $h\in \lie h_{aff}$ such that $h.v_\pi \neq 0$. Hence $\bold{0}\in G_\pi$. 
As the $\bold{A}_\Lambda$-module $A_\Lambda^\pi$ generated by $v_\pi$ corresponds to an algebra homomorphism on $\bold{A}_\Lambda$, the set $G_\pi$ is closed under addition. By \ref{highest.wt} every algebra homomorphism from $\bold{A}_\Lambda$ to $\C$ corresponds to an irreducible representation of $\bold{A}_\Lambda$. Hence $A_\Lambda^\pi$ is an irreducible $\bold{A}_\Lambda$-module and consequently every $\Z^{k-1}$-graded element of $A_\Lambda^\pi$ is invertible. That is if $A_\Lambda^\pi[\bom]\neq 0$, then $A_\Lambda^\pi [-\bom]\neq 0$. This implies that $G_\pi$ is closed under inverses and hence $G_\pi$ is a subgroup of $\Z^{k-1}$.

Suppose $M_1,M_2,\cdots,M_l$ is an enumeration of $\supp(\pi)$ then for any $h\in \lie h_{aff}$ and $2\leq i\leq k$, we have
$$ h\otimes t_i^r.v_\pi = (\sum_{i=1}^l \pi(M_i)(h) \ev_{M_i}(t_i^r)) v_\pi\otimes t_i^r.$$ Since $\supp(\pi)$ is finite and $\pi(M_i)(h)\in \Z_+$ for all $M_i\in \supp(\pi)$ and $h\in \lie h_{aff}$, $(\sum_{i=1}^l \pi(M_i)(h) \ev_{M_i}(t_i^r))$ cannot be zero for all $h\in \lie h_{aff}$ and $r\in \Z_+$. Hence for all $\pi\in \Pi$, rank of $G_\pi$ is $k-1$, implying that $G_\pi$ is a subgroup of $\Z^{k-1}$ of finite index.  \endproof
%Let $F_{k-1}$ be the set of all subgroups of $\Z^k-1$ of finite index and let $G: \Pi \rightarrow F_{k-1} $
Given $\pi\in \Pi$, we shall henceforth refer to the set $G_\pi$ as the group associated to $\pi.$
\vspace{.15cm}

Set $\C[t^{G_\pi}]$ as the set of all polynomials  in the variables $\{t^\bom: \bom\in G_\pi\}$. As $G_{\pi}$ is a subgroup of $\Z^{k-1}$ of rank $k-1$, it is easy to see that $\C[t^{G_\pi}]$ is isomorphic to a Laurent polynomial ring in $k-1$ variables. Let $G^\pi = \Z^{k-1}/G_\pi$. Clearly $G^\pi$ is a finite group.
For $\bog \in G^\pi$, let $\C t^\bog[t^{G_\pi}]$ be the set of all polynomials in the variables $\{t^{\bog+\bom} : \bom\in G_\pi\}$. From the construction of $G_\pi$ it follows that the irreducible $\bold{A}_\Lambda$-module $A_\Lambda^\pi$ is isomorphic to $v_\pi\otimes \C[t^{G_\pi}]$  and for each $\bog\in G^\pi$, the irreducible $\bold{A}_\Lambda$-module generated by the vector  
$v_\pi\otimes t^\bog \in L(X_\pi)$ is isomorphic to the subspace $v_\pi\otimes \C t^\bog[t^{G_\pi}]$ of $L(X_\pi)$.

The following result was proved in \cite[Proposition 3.5, Theorem 3.18, Example 4.2]{R2}.  

\begin{prop} For $\pi\in \Pi$, let $v_\pi$ be the highest weight vector of $X_\pi$ and let $G_\pi$ be the group associated to $\pi$ and $G^\pi=\Z^{k-1}/G_\pi$. Then we have the following.
\item[i.] For each $\bog\in G^\pi$,  the $\cal T(\lie g)$-module
$$X_\pi^\bog = \cal U(\cal T(\lie g)).v_\pi\otimes t^\bog,$$
is an irreducible $\cal T(\lie g)$-module.
\item[ii.] $L(X_\pi)$ is completely reducible as a $\cal T(\lie g)$-module. In fact as a $\cal T(\lie g)$-module $L(X_\pi)$ is isomorphic to the direct sum of the irreducible $\cal T(\lie g)$-modules $X_\pi^\bog, \bog\in G^\pi$, that is,
$$L(X_\pi)\cong_{\cal T(\lie g)} \underset{\bog\in G^\pi}{\bigoplus} X^\bog_\pi.$$
Further if $V$ is an irreducible $\cal T(\lie g)$-module in $\cal I_{fin}^{(m\boe_1)}, m>0$, then upto twisting by one-dimensional $\cal T(\lie g)$-modules $V$ is isomorphic to $X_\pi^\bog$ for some $\pi\in \Pi$ with $\wt(\pi)(\alpha_{n+1})=m$ and $\bog\in G_\pi$.
\end{prop}

\subsection{} Notice that for $\bob=(b_2,\cdots,b_k)\in (\C^\ast)^{k-1}$ the map $s_\bob:\C[t_2^{\pm1},\cdots,t_k^{\pm1}] \rightarrow \C[t_2^{\pm1},\cdots,t_k^{\pm1}]$ given by $t_i\mapsto b_it_i$, for $i=2,\cdots,k$  is an isomorphism. % Denote by $\bob.M$ the image of a maximal ideal $M$ of $\C[t_2^\pm,\cdots,t_k^\pm]$.
Using this isomorphism we now give a parametrization of the irreducible representations of $\cal T(\lie g)$ in $\cal I_\fn^{(m\boe_1)}, m>0$. In the case when the center acts trivially, the isomorphism classes of irreducible $\cal T(\lie g)$-modules in $\cal I_\fn$ was determined in  \cite{CFK,PB,L}.  

\label{isomorphism}

\begin{thm} Given $\pi,\pi'\in \Pi$, $\bog\in \Z^{k-1}/G_\pi$ and $\bog'\in \Z^{k-1}/G_{\pi'}$, the irreducible $\cal T(\lie g)$-module $X_\pi^\bog$ is isomorphic to $X_{\pi'}^{\bog'}$ if and only if  there exists $\bob=(b_2,\cdots,b_k)\in (\C^\ast)^{k-1}$ such that
\item[(i)]$\supp(\pi')= \{s_\bob(M): M\in \supp(\pi)\},$
\item[(ii)] For all $M\in \supp(\pi)$, %there exists one-dimensional $\lie g_{aff}$-module $Z_M$ such that  $X(\pi(M))\otimes Z_M$ is isomorphic to $X(\pi'(\bob.M))$ as a $\lie g_{aff}$-module or equivalently 
$X(\pi(M))$ is isomorphic to $X(\pi'(s_\bob(M)))$ as a $\lie g\otimes \C[t_1,t_1^{-1}]\oplus \C K_1$-module.  
 \item[(iii)]$\bog \equiv \bog'\mod G_\pi$.
\end{thm}
\proof Suppose $\chi:X^\bog_\pi \rightarrow X^{\bog'}_{\pi'}$ is a $\cal T(\lie g)$-module isomorphism from $X^\bog_\pi$ to $X^{\bog'}_{\pi'}$. If $M_1,M_2,\cdots, M_r$ is an enumeration of $\supp(\pi)$, then by \propref{Ch.Rmn.Thm} $X_\pi^\bog$ is a module for the Lie algebra $$ \underset{i=1}{\overset{r}{\oplus}} \big((\lie g_{\fn}\otimes \C[t_1^{\pm1}]) \otimes (\C[t_2^{\pm1},\cdots,t_k^{\pm1}]/M_i)\big) \oplus \cal Z\oplus D_k.$$ Hence via the isomorphism $\chi$, $X_{\pi'}^{\bog'}$ is a module for the Lie algebra  $\underset{i=1}{\overset{r}{\oplus}} \big((\lie g_{\fn}\otimes \C[t_1^{\pm1}]) \otimes (\C[t_2^{\pm1},\cdots,t_k^{\pm1}]/M_i)\big) \oplus \cal Z\oplus D_k.$ By definition however $X_{\pi'}^{\bog'}$ is a module for the Lie algebra  $${\underset{\cal M \in \supp(\pi')}{\oplus}} \big((\lie g_{\fn}\otimes \C[t_1^{\pm1}]) \otimes (\C[t_2^{\pm1},\cdots,t_k^{\pm1}]/\cal M)\big) \oplus \cal Z\oplus D_k.$$ Hence there exists an isomorphism $\sigma: \C[t_2^{\pm1},\cdots,t_k^{\pm1}]\rightarrow \C[t_2^{\pm1},\cdots,t_k^{\pm1}]$ such that for $1\leq i\leq r$, $\sigma(M_i) = \cal M_j^i$ for  $\cal M_j^i\in \supp(\pi'),$ or equivalently there exists $\bob=(b_2,\cdots,b_k)\in (\C^\times)^{k-1},$ such that  $$\sigma = s_\bob:\C[t_2^{\pm1},\cdots,t_k^{\pm1}]\rightarrow \C[t_2^{\pm1},\cdots,t_k^{\pm1}]$$ is the isomorphism of $\C[t_2^{\pm1},\cdots,t_k^{\pm1}]$ given by $t_i\mapsto b_it_i$ for $2\leq i\leq k. $ This implies $\{s_\bob(M): M\in \supp(\pi)\}\subseteq \supp(\pi').$ Since $\chi$ and $s_\bob$ are isomorphisms, taking the inverse maps we see that $$\supp(\pi') =\{s_\bob(M): M\in \supp(\pi)\}.$$

Let  $f_1,\cdots,f_r \in \C[t_2^{\pm1},\cdots,t_k^{\pm1}]$ be such that $f_i\in
\underset{j=1,j\neq i}{\overset{r}{\cap}} M_j$. Then for each
$x\in \lie g_\fn\otimes \C[t_1^{\pm1}]\oplus \C K_1 \oplus \C d_1$ we have,
$$ \chi(x\otimes f_i. v_\pi \otimes t^\bog) = x\otimes s_\bob(f_i).\chi(v_\pi \otimes t^\bog) \hspace{.5cm} \text{for}\ 1\leq i\leq r.$$ Using the definition of the action of $\cal T(\lie g)$ on $L(X_\pi)$ it follows that if $\chi$ is an isomorphism, then the $\lie g_{aff}$ module generated by $v_{M_i}$ is isomorphic to the $\lie g_{aff}$-module generated by $v_{s_\bob(M_i)}$.  Hence by
\cite[Theorem 3, Lemma 1]{VV} for every $M\in \supp(\pi)$, $X(\pi(M))$ is
isomorphic to $X(\pi'(s_\bob(M)))$ as a $\lie g_\fn\otimes \C[t_1,t_1^{-1}]\oplus \C K_1$-module.
 %there exists one-dimensional $\lie g_{aff}$ module $Z_{M}$ such that $X(\pi'(\bob.M))$ is isomorphic to $X(\pi(M))\otimes Z_M$ as a $\lie g_{aff}$ module.

As a consequence of conditions (i) and (ii), $G_\pi = G_{\pi'}$ whenever $X_\pi^\bog$ is isomorphic to $X_{\pi'}^{\bog'}$ implying that $\bog, \bog' \in \Z^{k-1}/G_\pi$. Since the set of highest weight vectors in $X_{\pi'}^{\bog'}$ are contained in the subspace $v_\pi'\otimes t^{\bog'}\C[t^{G_\pi'}] = v_\pi'\otimes t^{\bog'}\C[t^{G_\pi}]$, and the highest weight vector $v_\pi\otimes t^\bog$ of $X_\pi^\bog$ maps to a highest weight vector in $X_{\pi'}^{\bog'}$ we have, $ d_i(\chi(v_\pi\otimes t^\bog)) = \chi(d_i.v_\pi\otimes t^\bog) =
g_i\ \chi(v_\pi\otimes t^\bog), \hspace{.35cm}  \forall\ 2\leq i\leq k. $
This is possible if and only if $\bog \equiv \bog'\mod G_\pi$. This completes the proof of the theorem.
\endproof

\vspace{.25cm}

Using the same proof as above it is easy to see the following. 
\begin{prop}  Given $\pi,\pi'\in \Pi$, the irreducible $\lie g_\fn\otimes \C[t_1^{\pm1},\cdots,t_k^{\pm1}]\oplus \cal Z$-module $X_\pi$ is isomorphic to $X_{\pi'}$ if and only if  there exists $\bob=(b_2,\cdots,b_k)\in (\C^\ast)^{k-1}$ such that
\item[(i)]$\supp(\pi')= \{s_\bob(M): M\in \supp(\pi)\},$
\item[(ii)] For $M\in \supp(\pi)$, %there exists one-dimensional $\lie g_{aff}$-module $Z_M$ such that  $X(\pi(M))\otimes Z_M$ is isomorphic to $X(\pi'(\bob.M))$ as a $\lie g_{aff}$-module or equivalently 
$X(\pi(M))$ is isomorphic to $X(\pi'(\bob.M))$ as a $\lie g_\fn\otimes\C[t_1,t_1^{-1}]\oplus \C K_1$-modules.  
 \end{prop}

\end{document}